\documentclass[10pt]{article}
\usepackage{theorem}
\usepackage{amssymb}
\usepackage{latexsym}

\theorembodyfont{\itshape}
\newtheorem{Theorem}{Theorem}[section]
\newtheorem{s-Theorem}{Theorem}[subsection]
\newtheorem{Proposition}{Proposition}[section]
\newtheorem{s-Proposition}{Proposition}[subsection]

\newtheorem{s-Conjecture}{Conjecture}[subsection]
\newtheorem{Lemma}{Lemma}[section]
\newtheorem{s-Lemma}{Lemma}[subsection]

\newtheorem{s-Claim}{Claim}[subsection]

\newtheorem{Problem}{Problem}[section]
\newtheorem{Corollary}{Corollary}[section]
\newtheorem{s-Corollary}{Corollary}[subsection]

\theorembodyfont{\rmfamily}
\newtheorem{Definition}{Definition}[section]
\newtheorem{Notation}{Notation}[section]
\newtheorem{Remark}{Remark}[section]
\newtheorem{s-Remark}{Remark}[subsection]

\newtheorem{s-Example}{Example}[subsection]
\begin{document}
\title{{\bf A numerical characterization of polarized manifolds \textrm{\boldmath $(X,\mathcal{L})$} with \textrm{\boldmath $K_{X}=-(n-i)\mathcal{L}$} by the \boldmath{$i$}th sectional geometric genus and the \boldmath{$i$}th \boldmath{$\Delta$}-genus} 
\thanks{This research was partially supported by the Grant-in-Aid for Scientific Research (C) (No.20540045), Japan Society for the Promotion of Science, Japan.}}
\author{YOSHIAKI FUKUMA}
\date{}

\maketitle
\begin{abstract}
Let $(X,\mathcal{L})$ be a polarized manifold of dimension $n$.
In this paper, by using the $i$th sectional geometric genus and the $i$th $\Delta$-genus, we will give a numerical characterization of $(X,\mathcal{L})$ with $K_{X}=-(n-i)\mathcal{L}$ for the following cases (i) $i=2$, (ii) $i=3$ and $n\geq 5$, (iii) $\mbox{max}\{ 2, \dim \mbox{Bs}|\mathcal{L}|+2\}\leq i\leq n-1$.
\end{abstract}

\section{Introduction}
Let $X$ be a projective variety with $\dim X=n$ defined over 
the field of complex numbers, and let $\mathcal{L}$ be an ample line bundle on $X$.
Then $(X,\mathcal{L})$ is called a {\it polarized variety}.
If $X$ is smooth, then we say that $(X,\mathcal{L})$ is a polarized {\it manifold}. 
The main purpose of this paper is to give a numerical characterization of $(X,\mathcal{L})$ with $K_{X}=-(n-i)\mathcal{L}$.
Then the following is well-known:
\begin{Proposition}\label{PP0}
Let $(X,\mathcal{L})$ be a polarized manifold of dimension $n\geq 2$.
\begin{itemize}
\item [\rm (1)] $(X,\mathcal{L})$ is a polarized manifold with $K_{X}=-(n+1)\mathcal{L}$ {\rm (}resp. $K_{X}=-n\mathcal{L}${\rm )} if and only if $2g(X,\mathcal{L})-2=-2\mathcal{L}^{n}$ {\rm (}resp. $2g(X,\mathcal{L})-2=-\mathcal{L}^{n}${\rm )}.
\item [\rm (2)] {\rm (}See {\rm \cite[(1.9) Theorem]{Fujita80}}.{\rm )}
$(X,\mathcal{L})$ is a polarized manifold with $K_{X}=-(n-1)\mathcal{L}$, which is called a Del Pezzo manifold, if and only if $2g(X,\mathcal{L})-2=0$ and $\Delta(X,\mathcal{L})=1$. 
\end{itemize}
\end{Proposition}
\noindent
\\
(Here $g(X,\mathcal{L})$ (resp. $\Delta(X,\mathcal{L})$) denotes the sectional genus (resp. the $\Delta$-genus) of $(X,\mathcal{L})$.)
\par
As the next step, we want to give a numerical characterization of polarized manifolds with $K_{X}=-(n-i)\mathcal{L}$ and $i\geq 2$ by using some invariants of $(X,\mathcal{L})$.
In \cite{Fukuma04} and \cite{Fukuma05}, we define the $i$th sectional geometric genus $g_{i}(X,\mathcal{L})$ and the $i$th $\Delta$-genus $\Delta_{i}(X,\mathcal{L})$ of $(X,\mathcal{L})$ for every integer $i$ with $0\leq i\leq n$.
The $i$th sectional geometric genus (resp. the $i$th $\Delta$-genus) is a generalization of the sectional genus (resp. $\Delta$-genus), that is, $g_{1}(X,\mathcal{L})=g(X,\mathcal{L})$ (resp. $\Delta_{1}(X,\mathcal{L})=\Delta(X,\mathcal{L})$).
So by looking at the Proposition \ref{PP0} above carefully, the author thought maybe we were able to give a numerical characterization of polarized manifolds $(X,\mathcal{L})$ with $K_{X}=-(n-i)\mathcal{L}$ by using the $i$th sectional geometric genus and the $i$th $\Delta$-genus.
\par
In this paper, as the main results, we prove the following:
\begin{Theorem}\label{MainTheorems}
{\rm (}See Theorems {\rm \ref{MT4}}, {\rm \ref{MT2}} and {\rm \ref{MT3}} below.{\rm )}
Let $(X,\mathcal{L})$ be a polarized manifold of dimension $n\geq 3$.
Assume that one of the following types holds:
\begin{itemize}
\item [\rm (a)] $i=2$.
\item [\rm (b)] $i=3$ and $n\geq 5$.
\item [\rm (c)] $\mbox{\rm max}\{ 2, \dim\mbox{\rm Bs}|\mathcal{L}|+2 \}\leq i\leq n-1$.
\end{itemize}
Then the following are equivalent one another.
\begin{itemize}
\item [\rm $C(i,1)$] $K_{X}+(n-i)\mathcal{L}=\mathcal{O}_{X}$.
\item [\rm $C(i,2)$] $\Delta_{i}(X,\mathcal{L})=1$ and $2g_{1}(X,\mathcal{L})-2=(i-1)\mathcal{L}^{n}$.
\item [\rm $C(i,3)$] $\Delta_{i}(X,\mathcal{L})>0$ and $2g_{1}(X,\mathcal{L})-2=(i-1)\mathcal{L}^{n}$.
\item [\rm $C(i,4)$] $g_{i}(X,\mathcal{L})=1$ and $2g_{1}(X,\mathcal{L})-2=(i-1)\mathcal{L}^{n}$.
\item [\rm $C(i,5)$] $g_{i}(X,\mathcal{L})>0$ and $2g_{1}(X,\mathcal{L})-2=(i-1)\mathcal{L}^{n}$.
\end{itemize}
\end{Theorem}
\noindent
\par
The author would like to thank Dr. Hironobu Ishihara for giving some comments about this paper.

\begin{center}
{\bf Notation and Conventions}
\end{center}\noindent
\par
We say that $X$ is a {\it variety} if $X$ is an integral separated scheme 
of finite type.
In particular $X$ is irreducible and reduced if $X$ is a variety.
Varieties are always assumed to be defined 
over the field of complex numbers.
In this article, we shall study mainly a smooth projective variety.
The words ``line bundles" and ``Cartier divisors" are used interchangeably.
The tensor products of line bundles are denoted additively.
\\
$\mathcal{O}(D)$: invertible sheaf associated with a Cartier divisor $D$ on $X$.\\
$\mathcal{O}_{X}$: the structure sheaf of $X$.\\
$\chi(\mathcal{F})$:  the Euler-Poincar\'e characteristic of 
a coherent sheaf $\mathcal{F}$.\\
$h^{i}(\mathcal{F}):=\mbox{\rm dim}H^{i}(X,\mathcal{F})$ 
for a coherent sheaf $\mathcal{F}$ on $X$.
\\
$h^{i}(D):=h^{i}(\mathcal{O}(D))$ for a Cartier divisor $D$.
\\
$q(X)(=h^{1}(\mathcal{O}_{X}))$: the irregularity of $X$.
\\
$h^{i}(X,\mathbb{C}):=\dim H^{i}(X,\mathbb{C})$.
\\
$b_{i}(X):=h^{i}(X,\mathbb{C})$.
\\
$K_{X}$:  the canonical divisor of $X$.
\\
$\mathbb{P}^{n}$: the projective space of dimension $n$.
\\
$\mathbb{Q}^{n}$: a quadric hypersurface in $\mathbb{P}^{n+1}$.
\\
$\sim$ (or $=$):  linear equivalence.
\\
$\mbox{det}(\mathcal{E}):=\wedge^{r} \mathcal{E}$, 
where $\mathcal{E}$ is a vector bundle of rank $r$ on $X$.
\\
$\mathbb{P}_{X}(\mathcal{E})$: the projective space bundle 
associated with a vector bundle $\mathcal{E}$ on $X$.
\\
$H(\mathcal{E})$: 
the tautological line bundle on $\mathbb{P}_{X}(\mathcal{E})$.
\\
$\mathcal{F}^{\vee}:=Hom_{\mathcal{O}_{X}}(\mathcal{F},\mathcal{O}_{X})$.
\\
$c_{i}(\mathcal{E})$: the $i$-th Chern class of a vector bundle $\mathcal{E}$.
\\
$c_{i}(X):=c_{i}(\mathcal{T}_{X})$, where $\mathcal{T}_{X}$ is the tangent bundle of a smooth projective variety $X$.
\\
For a real number $m$ and a non-negative integer $n$, let 
\[
[m]^{n}:=\left\{
\begin{array}{cl}
m(m+1)\cdots (m+n-1) & \mbox{if $n\geq 1$,} \\
1 & \mbox{if $n=0$.}
\end{array} \right. \]
\[
[m]_{n}:=\left\{
\begin{array}{cl}
m(m-1)\cdots (m-n+1) & \mbox{if $n\geq 1$,} \\
1 & \mbox{if $n=0$.}
\end{array} \right. \]
Then for $n$ fixed, $[m]^{n}$ and $[m]_{n}$ are polynomials in $m$ whose degree are $n$.
\\
For any non-negative integer $n$,
\[
n!:=\left\{
\begin{array}{cl}
[n]_{n} & \mbox{if $n\geq 1$,} \\
1 & \mbox{if $n=0$.}
\end{array} \right. \]
Assume that $m$ and $n$ are integers with $n\geq 0$. Then we put
\[
{m\choose n}:=\frac{[m]_{n}}{n!} 
\]
We note that ${m\choose n}=0$ if $0\leq m<n$, and ${m\choose 0}=1$.

\section{Preliminaries}\label{S2}

Here we list up some facts which will be used later.

\begin{Definition}\label{DF3}
(1) Let $X$ (resp. $Y$) be an $n$-dimensional projective manifold, and let $\mathcal{L}$ (resp. $\mathcal{A}$) be an ample line bundle on $X$ (resp. $Y$).
Then $(X,\mathcal{L})$ is called a  {\it simple blowing up of $(Y,\mathcal{A})$} if there exists a birational morphism $\pi: X\to Y$ such that $\pi$ is a blowing up at a point of $Y$ and $\mathcal{L}=\pi^{*}(\mathcal{A})-E$, where $E$ is the exceptional divisor.
\\
(2) Let $X$ (resp. $M$) be an $n$-dimensional projective manifold, and let $\mathcal{L}$ (resp. $\mathcal{A}$) be an ample line bundle on $X$ (resp. $M$).
Then we say that $(M,\mathcal{A})$ is a {\it reduction of $(X,\mathcal{L})$} if $(X,\mathcal{L})$ is obtained by a composite of simple blowing ups of $(M,\mathcal{A})$, and $(M,\mathcal{A})$ is not obtained by a simple blowing up of any polarized manifold.
The morphism $\mu:X\to M$ is called the {\it reduction map}.
\end{Definition}

\begin{Definition}\label{DF4}
Let $(X,\mathcal{L})$ be a polarized manifold of dimension $n$.
We say that $(X,\mathcal{L})$ is a {\it scroll} (resp. {\it quadric fibration}, {\it Del Pezzo fibration}) {\it over a normal projective variety $Y$ with $\dim Y=m$} if there exists a surjective morphism with connected fibers $f:X\to Y$ such that $K_{X}+(n-m+1)\mathcal{L}=f^{*}\mathcal{A}$ (resp. $K_{X}+(n-m)\mathcal{L}=f^{*}\mathcal{A}$, $K_{X}+(n-m-1)\mathcal{L}=f^{*}\mathcal{A}$) for some ample line bundle $\mathcal{A}$ on $Y$.
\end{Definition}

\begin{Remark}\label{RM3}
If $(X,\mathcal{L})$ is a scroll over a smooth curve $C$ 
(resp. a smooth projective surface $S$) with $\dim X=n\geq 3$,
then by \cite[(3.2.1) Theorem]{BeSoWi} and \cite[Proposition 3.2.1 and Theorem 14.1.1]{BeSoBook}
there exists an ample vector bundle $\mathcal{E}$ of rank $n$ (resp. $n-1$)
on $C$ (resp. $S$) such that 
$(X,\mathcal{L})\cong (\mathbb{P}_{C}(\mathcal{E}),H(\mathcal{E}))$ 
(resp. $(\mathbb{P}_{S}(\mathcal{E}),H(\mathcal{E}))$).
\end{Remark}

\begin{Theorem}\label{TH1}
Let $(X,\mathcal{L})$ be a polarized manifold with $\dim X=n\geq 3$.
Then $(X,\mathcal{L})$ is one of the following types.
\begin{itemize}
\item [\rm (1)]
$(\mathbb{P}^{n}, \mathcal{O}_{\mathbb{P}^{n}}(1))$.
\item [\rm (2)]
$(\mathbb{Q}^{n}, \mathcal{O}_{\mathbb{Q}^{n}}(1))$.
\item [\rm (3)]
A scroll over a smooth projective curve.
\item [\rm (4)]
$K_{X}\sim -(n-1)\mathcal{L}$, that is, $(X,\mathcal{L})$ is a Del Pezzo manifold.
\item [\rm (5)]
A quadric fibration over a smooth curve.
\item [\rm (6)]
A scroll over a smooth projective surface.
\item [\rm (7)]
Let $(M, \mathcal{A})$ be a reduction of $(X,\mathcal{L})$.
\begin{itemize}
\item [{\rm (7.1)}] $n=4$, $(M, \mathcal{A})=(\mathbb{P}^{4}, \mathcal{O}_{\mathbb{P}^{4}}(2))$.
\item [{\rm (7.2)}]
$n=3$, $(M, \mathcal{A})=(\mathbb{Q}^{3}, \mathcal{O}_{\mathbb{Q}^{3}}(2))$.
\item [{\rm (7.3)}]
$n=3$, $(M, \mathcal{A})=(\mathbb{P}^{3}, \mathcal{O}_{\mathbb{P}^{3}}(3))$.
\item [{\rm (7.4)}]
$n=3$, $M$ is a $\mathbb{P}^{2}$-bundle over a smooth curve $C$ and 
for any fiber $F^{\prime}$ of it, $(F^{\prime}, \mathcal{A}|_{F^{\prime}})\cong(\mathbb{P}^{2}, \mathcal{O}_{\mathbb{P}^{2}}(2))$.
\item [\rm (7.5)]
$K_{M}\sim -(n-2)\mathcal{A}$, that is, 
$(M,\mathcal{A})$ is a Mukai manifold.
\item [\rm (7.6)]
$(M, \mathcal{A})$ is a Del Pezzo fibration over a smooth curve.
\item [\rm (7.7)]
$(M, \mathcal{A})$ is a quadric fibration over a normal surface.
\item [\rm (7.8)]
$n\geq 4$ and $(M,\mathcal{A})$ is a scroll over a normal projective variety of dimension $3$.
\item [\rm (7.9)] $K_{M}+(n-2)\mathcal{A}$ is nef and big.
\end{itemize}
\end{itemize}
\end{Theorem}
\noindent{\em Proof.}
See \cite[Proposition 7.2.2, Theorem 7.2.4, Theorem 7.3.2, Theorem 7.3.4, 
and Theorem 7.5.3]{BeSoBook}. 
See also \cite[Chapter II, (11.2), (11.7), and (11.8)]{FujitaBook}. $\Box$

\begin{Remark}\label{RM2}
Let $(X,\mathcal{L})$ be a polarized manifold with $\dim X=n\geq 3$.
\begin{itemize}
\item [(1)]
$\kappa(K_{X}+(n-2)\mathcal{L})=-\infty$ if and only if $(X,\mathcal{L})$ is one of the types from (1) to (7.4) in Theorem \ref{TH1}.
\item [(2)]
$\kappa(K_{X}+(n-2)\mathcal{L})=0$ if and only if $(X,\mathcal{L})$ is (7.5) in Theorem \ref{TH1}.
\item [(3)]
$\kappa(K_{X}+(n-2)\mathcal{L})\geq 1$ if and only if $(X,\mathcal{L})$ is one of the types from (7.6) to (7.9) in Theorem \ref{TH1}.
\end{itemize}
\end{Remark}

\begin{Definition}\label{DF5}(\cite[7.5.7 Definition-Notation]{BeSoBook})
Let $(X,\mathcal{L})$ be a polarized manifold of dimension $n\geq 3$, and let $(M,\mathcal{A})$ be a reduction of $(X,\mathcal{L})$.
Assume that $K_{M}+(n-2)\mathcal{A}$ is nef and big.
Then for large $m\gg 0$ the morphism $\varphi: M\to W$ associated to $|m(K_{M}+(n-2)\mathcal{A})|$ has connected fibers and normal image $W$.
Then we note that there exists an ample line bundle $\mathcal{K}$ on $W$ such that $K_{M}+(n-2)\mathcal{A}=\varphi^{*}(\mathcal{K})$.
Let $\mathcal{D}:=(\varphi_{*}\mathcal{A})^{\vee\vee}$, where $^{\vee\vee}$ denotes the double dual.
Then the pair $(W,\mathcal{D})$ together with $\varphi$ is called the {\it second reduction of $(X,\mathcal{L})$}. 
\end{Definition}

\begin{Remark}\label{RM6}
(1) If $K_{M}+(n-2)\mathcal{A}$ is nef and big but not ample, then $\varphi$ is equal to the nef value morphism of $\mathcal{A}$.
\\
(2) If $K_{M}+(n-2)\mathcal{A}$ is ample, then $\varphi$ is an isomorphism.
\\
(3) If $n\geq 4$, then $W$ has isolated terminal singularities and is $2$-factorial.
Moreover if $n$ is even, then $X$ is Gorenstein (see \cite[Proposition 7.5.6]{BeSoBook}).
\end{Remark}

Here we consider the characterization of $(X,\mathcal{L})$ with $\kappa(K_{X}+(n-3)\mathcal{L})=-\infty$.
We note that $\kappa(K_{X}+(n-1)\mathcal{L})=-\infty$ (resp. $\kappa(K_{X}+(n-2)\mathcal{L})=-\infty$) if and only if $(X,\mathcal{L})$ is one of the types from (1) to (3) (resp. from (1) to (7.4)) in Theorem \ref{TH1}.
Here we consider the case where $\kappa(K_{X}+(n-3)\mathcal{L})=-\infty$.
If $(X,\mathcal{L})$ is one of the types from (1) to (7.8) in Theorem \ref{TH1}, then $\kappa(K_{X}+(n-3)\mathcal{L})=-\infty$ holds. 
So we assume that $K_{M}+(n-2)\mathcal{A}$ is nef and big.
Then there exist a normal projective variety $W$ with only 2-factorial isolated terminal singularities, a birational morphism $\phi_{2}: M\to W$ and an ample line bundle $\mathcal{K}$ on $W$ such that $K_{M}+(n-2)\mathcal{A}=(\phi_{2})^{*}(\mathcal{K})$.
Let $\mathcal{D}:=(\phi_{2})_{*}(\mathcal{A})^{\vee\vee}$.
Then $\mathcal{D}$ is a 2-Cartier divisor on $W$ and $\mathcal{K}=K_{W}+(n-2)\mathcal{D}$ (see \cite[Lemma 7.5.8]{BeSoBook}).
Then the pair $(W,\mathcal{D})$ is the second reduction of $(X,\mathcal{L})$
(see Definition \ref{DF5}).
Here we note that if $K_{M}+(n-2)\mathcal{A}$ is ample, then $(W,\mathcal{K})\cong (M, K_{M}+(n-2)\mathcal{A})$.
\par
Then the following properties hold:
\begin{itemize}
\item [(1)] $\kappa(K_{X}+(n-3)\mathcal{L})=\kappa(K_{W}+(n-3)\mathcal{K})$ holds \cite[Corollary 7.6.2]{BeSoBook}.
\item [(2)] $(n-2)(K_{W}+(n-3)\mathcal{D})=K_{W}+(n-3)\mathcal{K}$ and $K_{M}+(n-3)\mathcal{A}=\phi_{2}^{*}(K_{W}+(n-3)\mathcal{D})+\Delta$ for an exceptional $\mathbb{Q}$-effective divisor $\Delta$ of $\phi_{2}$.
Therefore
\begin{eqnarray*}
m(n-2)(K_{X}+(n-3)\mathcal{L})
&=&m(n-2)\phi_{1}^{*}(K_{M}+(n-3)\mathcal{A})+E_{1} \\
&=&m(n-2)\phi_{1}^{*}\circ \phi_{2}^{*}(K_{W}+(n-3)\mathcal{D})+E_{1}+m(n-2)\Delta\\
&=&m\phi_{1}^{*}\circ \phi_{2}^{*}(K_{W}+(n-3)\mathcal{K})+E_{1}+m(n-2)\Delta.
\end{eqnarray*}
(Here $\phi_{1}:X\to M$ is a reduction of $(X,\mathcal{L})$ and $E_{1}$ is a $\phi_{1}$-exceptional effective divisor.)
\item [(3)] $h^{0}((n-2)m(K_{X}+(n-3)\mathcal{L}))=h^{0}(m(K_{W}+(n-3)\mathcal{K}))$ for every integer $m$ with $m\geq 1$.
\end{itemize}

Moreover if $n\geq 4$, then there exists a normal factorial projective variety $M^{\sharp}$ with only isolated terminal singularities and birational morphisms $\phi_{2}^{\sharp}: M\to M^{\sharp}$ and $\psi: M^{\sharp}\to W$ such that $\phi_{2}=\psi\circ\phi_{2}^{\sharp}$.
Then $M^{\sharp}$ is called the {\it factorial stage} (see \cite[7.5.7 Definition-Notation]{BeSoBook} or \cite[(2.6)]{Fujita92}).
\par
If $\tau(\mathcal{K})\leq n-3$, then by above we see that $\kappa(K_{X}+(n-3)\mathcal{L})\geq 0$. (Here $\tau(\mathcal{K})$ denotes the nef value of $\mathcal{K}$.)
So we may assume that $\tau(\mathcal{K})>n-3$.
\par
If $n\geq 5$, then $(W, \mathcal{K})$ with $\tau(\mathcal{K})>n-3$ is one of some special types by \cite[Theorems 7.7.2, 7.7.3 and 7.7.5]{BeSoBook}.
So we can get the following:
\begin{Proposition}\label{PP2}
Let $(X,\mathcal{L})$ be a polarized manifold of dimension $n\geq 5$, $(M,\mathcal{A})$ a reduction of $(X,\mathcal{L})$, and $(W,\mathcal{K})$ the second reduction of $(X,\mathcal{L})$.
Then $\kappa(K_{X}+(n-3)\mathcal{L})=-\infty$ if and only if $(X,\mathcal{L})$ satisfies one of the following:
\begin{itemize}
\item [\rm (1)] $(X,\mathcal{L})$ is one of the types {\rm (1), (2), (3), (4), (5), (6), (7.5), (7.6), (7.7)} or {\rm (7.8)} in Theorem {\rm \ref{TH1}}.
\item [\rm (2)] $K_{M}+(n-2)\mathcal{A}$ is nef and big, and $(W,\mathcal{K})$ is one of the following:
\begin{itemize}
\item [\rm (2.1)] $(\mathbb{P}^{6}, \mathcal{O}_{\mathbb{P}^{6}}(1))$.
\item [\rm (2.2)] $1$, $2$ or $3$ in {\rm \cite[Theorem {\rm 7.7.5}]{BeSoBook}}.
\end{itemize}
\end{itemize}
\end{Proposition}
\noindent{\em Proof.}
See \cite[Theorems 7.7.2, 7.7.3, 7.7.5 and Proposition 7.7.9]{BeSoBook}. $\Box$
\\
\par
So we consider the case of $n=4$.
In this case $M^{\sharp}$ and $W$ are Gorenstein (see \cite[Proposition 7.5.6 and 7.5.7 Definition-Notation]{BeSoBook}).
Then by the proof of \cite[Section 4]{Fujita92} we see that $(W,\mathcal{K})$ or $M^{\sharp}$ is one of the types in \cite[(4.$\infty$)]{Fujita92}.
If $(W,\mathcal{K})$ or $M^{\sharp}$ is either (4.2), (4.4.0), (4.4.1), (4.4.2), (4.6.0.0), (4.6.0.1.0), (4.6.0.2.1), (4.6.1), (4.7) or (4.8.0) in \cite[(4.$\infty$)]{Fujita92}, then we see that $\kappa(K_{X}+\mathcal{L})=-\infty$.\par
Assume that $(W,\mathcal{K})$ is the type (4.4.4) in \cite[(4.$\infty$)]{Fujita92}. Then we note that $\tau(\mathcal{K})=3$ and there exist a normal Gorenstein projective variety $W_{2}$, an ample line bundle $\mathcal{K}_{2}$ on $W_{2}$ and a birational morphism $\mu: W\to W_{2}$ such that $\mu$ is the simultaneous contraction to distinct smooth points of divisors $E_{i}\cong \mathbb{P}^{3}$ such that $E_{i}\subset \mbox{reg}(W)$, $E_{i}|_{E_{i}}\cong\mathcal{O}_{\mathbb{P}^{3}}(-1)$, $K_{W}+3\mathcal{K}=\mu^{*}(K_{W_{2}}+3\mathcal{K}_{2})$ and $K_{W_{2}}+3\mathcal{K}_{2}$ is ample, that is, $\tau(\mathcal{K}_{2})<3$.
Moreover we infer that $W_{2}$ has the same singularities as $W$ by above.
Since $E_{i}\subset \mbox{reg}(W)$, we have $\psi^{-1}(E_{i})\cong E_{i}$ by the definition of $\psi$.
Hence there exist a normal Gorenstein projective variety $W_{2}^{\sharp}$ and birational morphisms $\mu^{\sharp}: M^{\sharp}\to W_{2}^{\sharp}$ and $\psi^{\sharp}:W_{2}^{\sharp}\to W_{2}$ such that $\mu\circ\psi=\psi^{\sharp}\circ\mu^{\sharp}$.
We note that $\mu^{\sharp}: M^{\sharp}\to W_{2}^{\sharp}$ is the contraction of $\psi^{-1}(E_{i})$ and $W_{2}^{\sharp}$ has the same singularities as $M^{\sharp}$.
The pair $(W_{2},\mathcal{K}_{2})$ is a reduction of $(W,\mathcal{K})$ and is called the {\it $2\frac{1}{2}$ reduction of $(W,\mathcal{K})$} in \cite[(2.2) Theorem-Definition]{BeSo93}.
We also note that $h^{j}(\mathcal{O}_{X})=h^{j}(\mathcal{O}_{M})=h^{j}(\mathcal{O}_{W})=h^{j}(\mathcal{O}_{W_{2}})=h^{j}(\mathcal{O}_{M^{\sharp}})=h^{j}(\mathcal{O}_{W^{\sharp}})$.
For this $\psi^{\sharp}:W_{2}^{\sharp}\to W_{2}$ and $(W_{2},\mathcal{K}_{2})$, we can apply the same argument as in \cite[Section 4]{Fujita92}.
If $\tau(\mathcal{K}_{2})\leq 1$, then we can prove that $\kappa(K_{X}+\mathcal{L})\geq 0$.
So we assume that $\tau(\mathcal{K}_{2})>1$.
Then $(W_{2},\mathcal{K}_{2})$ is either (4.6.0.0), (4.6.0.1.0), (4.6.0.2.1), (4.6.1), (4.6.4), (4.7) or (4.8.0) in \cite[(4.$\infty$)]{Fujita92}.
\par
If $(W_{2},\mathcal{K}_{2})$ is either (4.6.0.0), (4.6.0.1.0), (4.6.0.2.1), (4.6.1), (4.7) or (4.8.0) in \cite[(4.$\infty$)]{Fujita92}, then we see that $\kappa(K_{X}+\mathcal{L})=-\infty$.
\par
If $(W_{2},\mathcal{K}_{2})$ is the type (4.6.4) in \cite[(4.$\infty$)]{Fujita92}, then by the same argument as in \cite[Section 4]{Fujita92} we see that there exist a normal Gorenstein projective variety $W_{3}$, an ample line bundle $\mathcal{K}_{3}$ on $W_{3}$ and a birational morphism $\mu_{2}: W_{2}\to W_{3}$ such that $W_{3}$ has the same singularities as $W_{2}$, $K_{W_{2}}+2\mathcal{K}_{2}=\mu_{2}^{*}(K_{W_{3}}+2\mathcal{K}_{3})$ and 
$K_{W_{3}}+2\mathcal{K}_{3}$ is ample, that is, $\tau(\mathcal{K}_{3})<2$.
Here we note that $\kappa(K_{X}+\mathcal{L})=\kappa(K_{W_{2}}+\mathcal{K}_{2})=\kappa(K_{W_{3}}+\mathcal{K}_{3})$.
\par
If $\tau(\mathcal{K}_{3})\leq 1$, 
then $\kappa(K_{X}+\mathcal{L})=\kappa(K_{W_{3}}+\mathcal{K}_{3})\geq 0$.
\par
If $\tau(\mathcal{K}_{3})>1$, then $(W_{3},\mathcal{K}_{3})$ is either (4.7) or (4.8.0) by the same argument as in \cite[Section 4]{Fujita92} and we have $\kappa(K_{X}+\mathcal{L})=\kappa(K_{W_{3}}+\mathcal{K}_{3})=-\infty$.
\\
\par
By the above argument, we get the following:
\begin{Theorem}\label{TH2}
Let $(X,\mathcal{L})$ be a polarized manifold of dimension $n=4$.
\\
\\
{\rm (1)} The inequality $\kappa(K_{X}+\mathcal{L})\geq 0$ holds if and only if there exist a normal projective variety $W_{3}$ with only isolated terminal singularities, an ample line bundle $\mathcal{K}_{3}$ on $W_{3}$, and a birational morphism $\Phi: X\to W_{3}$ such that $\tau(\mathcal{K}_{3})\leq 1$ and $h^{0}(2m(K_{X}+\mathcal{L}))=h^{0}(m(K_{W_{3}}+\mathcal{K}_{3}))$ for every positive integer $m$.
\\
\\
{\rm (2)} The equality $\kappa(K_{X}+\mathcal{L})=-\infty$ holds if and only if $(X,\mathcal{L})$ satisfies one of the following:
\begin{itemize}
\item [\rm (2.1)] The pair $(X,\mathcal{L})$ is either {\rm (1), (2), (3), (4), (5), (6), (7.1), (7.5), (7.6), (7.7)} or {\rm (7.8)} in Theorem {\rm \ref{TH1}}.
\item [\rm (2.2)] There exist a normal projective variety $W_{3}$ with only isolated terminal singularities, an ample line bundle $\mathcal{K}_{3}$ on $W_{3}$, and a birational morphism $\Phi: X\to W_{3}$ such that $(W_{3}, \mathcal{K}_{3})$ is either {\rm (4.2), (4.4.0), (4.4.1), (4.4.2), (4.6.0.0), (4.6.0.1.0), (4.6.0.2.1), (4.6.1), (4.7)} or {\rm (4.8.0)} in {\rm \cite[(4.$\infty$)]{Fujita92}}.
\end{itemize}
\end{Theorem}

Furthermore we need the following two lemmas.

\begin{Lemma}\label{LM1}
Let $(X,\mathcal{L})$ be a polarized manifold of dimension $n\geq 4$, and let $(M,\mathcal{A})$ be a reduction of $(X,\mathcal{L})$.
If $\kappa(K_{X}+(n-3)\mathcal{L})=-\infty$, then $h^{j}(\mathcal{O}_{X})=0$ for any $j$ with $j\geq 3$ unless $(M,\mathcal{A})$ is a scroll over a normal projective variety of dimension $3$.
If $(M,\mathcal{A})$ is a scroll over a normal projective variety of dimension $3$, then $h^{j}(\mathcal{O}_{X})=h^{j}(\mathcal{O}_{M})=0$ for every integer $j$ with $j\geq 4$.
\end{Lemma}
{\em Proof.}
(I) First we assume that $n\geq 5$.
By assumption and Proposition \ref{PP2}, $(X,\mathcal{L})$ satisfies either (1), (2.1) or (2.2) in Proposition \ref{PP2}.
Here we note that since $h^{j}(\mathcal{O}_{X})=h^{j}(\mathcal{O}_{M})=h^{j}(\mathcal{O}_{W})$, we have ony to prove that $h^{j}(\mathcal{O}_{W})=0$.
But by Proposition \ref{PP2} it is easy to show this and left to the reader.
\\
\\
(II) Next we assume that $n=4$.
By Theorem \ref{TH2}, $(X,\mathcal{L})$ satisfies either (2.1) or (2.2) in Theorem \ref{TH2}.
Here we use notation above.
Here we note that $h^{j}(\mathcal{O}_{X})=h^{j}(\mathcal{O}_{M})=h^{j}(\mathcal{O}_{W})=h^{j}(\mathcal{O}_{W_{2}})=h^{j}(\mathcal{O}_{W_{3}})=h^{j}(\mathcal{O}_{M^{\sharp}})=h^{j}(\mathcal{O}_{W^{\sharp}})$.
\\
\\
(II.A) If $(W_{3},\mathcal{K}_{3})$ is a $(\mathbb{P}^{3},\mathcal{O}_{\mathbb{P}^{3}}(2))$-fibration over a smooth curve $C$, then there exists an extremal ray $R$ such that $(K_{W_{3}}+2\mathcal{K}_{3})R=0$ (see \cite[(4.6.1)]{Fujita92}).
Let $\rho$ be the contraction morphism of $R$.
Then $\dim \rho(X)=1$, $\rho(X)=C$ and $\rho: X\to C$ is the $(\mathbb{P}^{3},\mathcal{O}_{\mathbb{P}^{3}}(2))$-fibration.
Moreover there exists a line bundle $\mathcal{B}$ on $C$ such that $K_{W_{3}}+2\mathcal{K}_{3}=\rho^{*}(\mathcal{B})$.
Since $\mathcal{K}_{3}$ is ample, by \cite[Theorem 1-2-5]{KaMaMa} we see that $R^{j}\rho_{*}(K_{W_{3}}+2\mathcal{K}_{3})=0$ for every integer $j$ with $j>0$.
Hence $0=R^{j}\rho_{*}(\rho^{*}(\mathcal{B}))\cong R^{j}\rho_{*}(\mathcal{O}_{W_{3}})\otimes \mathcal{B}$ and we have $R^{j}\rho_{*}(\mathcal{O}_{W_{3}})=0$ for every positive integer $j$.
Therefore $h^{j}(\mathcal{O}_{W_{3}})=h^{j}(\rho_{*}(\mathcal{O}_{W_{3}}))$.
On the other hand $h^{j}(\rho_{*}(\mathcal{O}_{W_{3}}))=0$ for every integer $j$ with $j\geq 2$.
Therefore $h^{j}(\mathcal{O}_{W_{3}})=0$ for every $j\geq 2$.
\\
\\
(II.B) Assume that $(W_{3},\mathcal{K}_{3})$ is the cone over $(\mathbb{P}^{3},\mathcal{O}_{\mathbb{P}^{3}}(2))$.
Then the following holds: Let $\mathcal{E}:=\mathcal{O}_{\mathbb{P}^{3}}\oplus \mathcal{O}_{\mathbb{P}^{3}}(2)$, $P:=\mathbb{P}_{\mathbb{P}^{3}}(\mathcal{E})$ and $H(\mathcal{E})$ the tautological line bundle on $P$.
Let $\pi: P\to \mathbb{P}^{N}$ be the morphism associated to $H(\mathcal{E})$.
Then $W_{3}=\pi(P)$ (see \cite[1.1.8 in Chapter I]{BeSoBook}).
First we note that $h^{j}(\mathcal{O}_{P})\geq h^{j}(\mathcal{O}_{W_{3}})$ for every nonnegative integer $j$.
On the other hand $h^{j}(\mathcal{O}_{P})=h^{j}(\mathcal{O}_{\mathbb{P}^{3}})=0$ for every integer $j$ with $j\geq 1$.
Therefore we get $h^{j}(\mathcal{O}_{W_{3}})=0$ for every integer $j$ with $j\geq 1$.
\\
\\
(II.C) For other cases it is easy and left to the reader. $\Box$

\begin{Lemma}\label{LM2}
Let $(X,\mathcal{L})$ be a polarized manifold of dimension $n\geq 3$, and let $(M,\mathcal{A})$ be a reduction of $(X,\mathcal{L})$.
Assume that $K_{M}+(n-2)\mathcal{A}$ is nef and big.
Let $(W,\mathcal{D})$ be the second reduction of $(X,\mathcal{L})$ and $\varphi:M\to W$ its morphism.
{\rm (}Here we use notation in Definition {\rm \ref{DF5}}.{\rm )}
Then $h^{j}(\mathcal{A})=h^{j}(\mathcal{D})$ for every integer $j\geq 3$.
\end{Lemma}
{\em Proof.}
As we said in Definition \ref{DF5}, there exists an ample line bundle $\mathcal{K}$ on $Y$ such that $K_{M}+(n-2)\mathcal{A}=\varphi^{*}(\mathcal{K})$.
By \cite[Theorem 1-2-5]{KaMaMa} we have $R^{j}\varphi_{*}(K_{M}+(n-1)\mathcal{A})=0$.
On the other hand $R^{j}\varphi_{*}\mathcal{O}(K_{M}+(n-1)\mathcal{A})=R^{j}\varphi_{*}(\varphi^{*}(\mathcal{K})\otimes \mathcal{A})=\mathcal{K}\otimes R^{j}\varphi_{*}(\mathcal{A})$.
Therefore $R^{j}\varphi_{*}(\mathcal{A})=0$ and we have $h^{j}(\mathcal{A})=h^{j}(\varphi_{*}(\mathcal{A}))$ for every positive integer $j$.
\par
Since $\mathcal{A}$ is a line bundle on $M$, we see that $\varphi_{*}(\mathcal{A})$ is a torsion free coherent sheaf on $W$.
Then there exists an injective homomorphism $\mu: \varphi_{*}(\mathcal{A})\to (\varphi_{*}(\mathcal{A}))^{\vee\vee}$.
Hence we get the following exact sequence
$$0\to \varphi_{*}(\mathcal{A})\to (\varphi_{*}(\mathcal{A}))^{\vee\vee}\to \mbox{Coker}\mu\to 0.$$
Here we note that $\dim \mbox{Supp}(\mbox{Coker}\mu)\leq 1$ because there exists a closed subset $Z$ on $W$ such that $\dim Z\leq 1$ and $M\backslash\varphi^{-1}(Z)\cong W\backslash Z$.
Therefore $h^{j}(\mbox{Coker}\mu)=0$ for every integer $j$ with $j\geq 2$ by \cite[Theorem 2.7 in Chapter III]{Hartshorne} or \cite[Theorem 4.6$^{*}$]{Iitaka}.
Hence we have $h^{j}(\varphi_{*}(\mathcal{A}))=h^{j}((\varphi_{*}(\mathcal{A}))^{\vee\vee})$ for every integer $j$ with $j\geq 3$.
Since $\mathcal{D}=(\varphi_{*}(\mathcal{A}))^{\vee\vee}$, we get the assertion. $\Box$

\begin{Definition}
Let $X$ be a smooth projective variety and let $\mathcal{F}$ be a vector bundle on $X$.
Then for every integer $j$ with $j\geq 0$, the {\it $j$-th Segre class $s_{j}(\mathcal{F})$ of $\mathcal{F}$} is defined by the following equation: $c_{t}(\mathcal{F}^{\vee})s_{t}(\mathcal{F})=1$, where $c_{t}(\mathcal{F}^{\vee})$ is the Chern polynomial of $\mathcal{F}^{\vee}$ and
$s_{t}(\mathcal{F})=\sum_{j\geq 0}s_{j}(\mathcal{F})t^{j}$.
\end{Definition}

\begin{Remark}
(1) Let $X$ be a smooth projective variety and let $\mathcal{F}$ be a vector bundle on $X$.
Let $\tilde{s}_{j}(\mathcal{F})$ be the Segre class which is defined in \cite[Chapter 3]{Fulton}.
Then $s_{j}(\mathcal{F})=\tilde{s}_{j}(\mathcal{F}^{\vee})$.
\\
(2) For every integer $i$ with $1\leq i$, $s_{i}(\mathcal{F})$ can be written by using the Chern classes $c_{j}(\mathcal{F})$ with $1\leq j\leq i$.
(For example, $s_{1}(\mathcal{F})=c_{1}(\mathcal{F})$, $s_{2}(\mathcal{F})=c_{1}(\mathcal{F})^{2}-c_{2}(\mathcal{F})$, and so on.)
\end{Remark}

\section{Review on the \boldmath{$i$}th sectional geometric genus and the \boldmath{$i$}th $\Delta$-genus of polarized varieties.}\label{S3}

Here we are going to review the $i$th sectional geometric genus and the $i$th $\Delta$-genus of polarized varieties $(X,\mathcal{L})$ for every integer $i$ with $0\leq i\leq \dim X$.
Up to now, there are many investigations of $(X,\mathcal{L})$ via the sectional genus and the $\Delta$-genus.
In order to analyze $(X,\mathcal{L})$ more deeply, the author extended these notions.
In \cite[Definition 2.1]{Fukuma04} we defined an invariant called the {\it $i$th sectional geometric genus} which is thought to be a generalization 
of the sectional genus.
First we recall the definition of this invariant.

\begin{Notation}\label{1.1}
Let $(X,\mathcal{L})$ be a polarized variety of dimension $n$, and let $\chi(t\mathcal{L})$ be the Euler-Poincar\'e characteristic of $t\mathcal{L}$.
Then $\chi(t\mathcal{L})$ is a polynomial in $t$ of degree $n$, and we can describe this as 
$$\chi(t\mathcal{L})=\sum_{j=0}^{n}\chi_{j}(X,\mathcal{L}){t+j-1\choose j}.$$
\end{Notation}

\begin{Definition}\label{1.2}
(\cite[Definition 2.1]{Fukuma04})
Let $(X,\mathcal{L})$ be a polarized variety of dimension $n$.
Then for any integer $i$ with $0\leq i\leq n$ the {\it $i$th sectional geometric genus $g_{i}(X,\mathcal{L})$ of $(X,\mathcal{L})$} is defined by the following.
$$g_{i}(X,\mathcal{L})=(-1)^{i}(\chi_{n-i}(X,\mathcal{L})-\chi(\mathcal{O}_{X}))+\sum_{j=0}^{n-i}(-1)^{n-i-j}h^{n-j}(\mathcal{O}_{X}).$$
\end{Definition}

\begin{Remark}\label{Remark1} 
\begin{itemize}
\item [(1)] 
Since $\chi_{n-i}(X,\mathcal{L})\in\mathbb{Z}$, the invariant $g_{i}(X,\mathcal{L})$ is an integer by definition.
\item [(2)] If $i=\dim X=n$, then $g_{n}(X,\mathcal{L})=h^{n}(\mathcal{O}_{X})$ and $g_{n}(X,\mathcal{L})$ is independent of $\mathcal{L}$.
\item [(3)] If $i=0$, then $g_{0}(X,\mathcal{L})=\mathcal{L}^{n}$.
\item [(4)] If $i=1$, then $g_{1}(X,\mathcal{L})=g(\mathcal{L})$, where $g(\mathcal{L})$ is the sectional genus of $(X,\mathcal{L})$.
If $X$ is smooth, then $g_{1}(X,\mathcal{L})=1+(1/2)(K_{X}+(n-1)\mathcal{L})\mathcal{L}^{n-1}$, where $K_{X}$ denotes the canonical line bundle on $X$.
\item [(5)] 
Let $(X,\mathcal{L})$ be a polarized manifold of dimension $n$ and let $(M,\mathcal{A})$ be a reduction of $(X,\mathcal{L})$.
Then $g_{i}(X,\mathcal{L})=g_{i}(M,\mathcal{A})$ for every integer $i$ with $1\leq i\leq n$.
\end{itemize}
\end{Remark}

The following are main problems about the $i$th sectional geometric genus.

\begin{Problem}\label{SGG-PB1}
\begin{itemize}
\item [\rm (i)] Does the $i$th sectional geometric genus have a property similar to that of the sectional genus $?$ For example, there are the following two questions.
\begin{itemize}
\item [\rm (i.1)] Does $g_{i}(X,\mathcal{L})\geq 0$ hold $?$
More strongly, does $g_{i}(X,\mathcal{L})\geq h^{i}(\mathcal{O}_{X})$ hold $?$
\item [\rm (i.2)] Can we get the $i$th sectional geometric genus version of the theory on sectional genus $?$
\end{itemize}
\item [\rm (ii)] Are there any relationship 
between $g_{i}(X,\mathcal{L})$ and $g_{i+1}(X,\mathcal{L})$ $?$
\item [\rm (iii)] Classify $(X,\mathcal{L})$ by the value of the $i$th sectional geometric genus.
\item [\rm (iv)] What is the geometric meaning of the $i$th sectional geometric genus $?$
\end{itemize}
\end{Problem}

\begin{Remark}
(1) First we consider Problem \ref{SGG-PB1} (i.1).
At present we can prove that the non-negativity of $g_{i}(X,\mathcal{L})$ holds if (a) $i=0$, (b) $i=1$, (c) $i=2$ and $n=3$, (d) $i=n$.
But in general it is unknown whether $g_{i}(X,\mathcal{L})$ is non-negative or not.
Next we consider the second inequality. Of course, if the second inequality holds, then non-negativity of $g_{i}(X,\mathcal{L})$ also holds.
If $i=0$ or $n$, then this inequality holds.
But it is unknown whether this inequality holds or not in general.
If $i=1$, then this is a conjecture proposed by Fujita \cite[(13.7) Remark]{FujitaBook}, and this case has been studied for several cases (see, for example, \cite{Fukuma97-1}, \cite{Fukuma97-2}, \cite{Fukuma99} and so on).
In \cite[Corollary 2.8]{Fukuma04-2}, we showed that the second inequality holds if $\dim\mbox{Bs}|\mathcal{L}|+1\leq i\leq n-1$.
\\
(2) Next we consider Problem \ref{SGG-PB1} (ii).
If $\mbox{Bs}|L|=\emptyset$, then $g_{i}(X,\mathcal{L})=0$ implies $g_{i+1}(X,\mathcal{L})=0$.
\\
(3) Next we consider Problem \ref{SGG-PB1} (iii).
If $i=1$, then the classification of polarized manifolds $(X,\mathcal{L})$ with $g_{1}(X,\mathcal{L})\leq 2$ was obtained (see \cite{Fujita87-1}, \cite{Ionescu86-2}, \cite{BeLaPa87}, and \cite{Fujita87-2}).
\par
If $i=2$, then the classification of polarized manifolds $(X,\mathcal{L})$ with the following is obtained (see \cite[Corollary 3.5 and Theorem 3.6]{Fukuma04} and \cite{Fukuma08-4}):
\begin{itemize}
\item [\rm (ii.1)] The case where $\mbox{Bs}|\mathcal{L}|=\emptyset$ and $g_{2}(X,\mathcal{L})=h^{2}(\mathcal{O}_{X})$.
\item [\rm (ii.2)] The case where $\mathcal{L}$ is very ample and $g_{2}(X,\mathcal{L})=h^{2}(\mathcal{O}_{X})+1$.
\end{itemize}
\noindent
\\
(4) Finally we consider Problem \ref{SGG-PB1} (iv).
Namely we will explain the geometric meaning of the $i$th sectional geometric genus.
First we will give the following definition.

\begin{Definition}\label{1.5}
Let $(X,\mathcal{L})$ be a polarized variety of dimension $n$.
Then $\mathcal{L}$ has a {\it $k$-ladder} if there exists a sequence of irreducible and reduced subvarieties $X\supset X_{1}\supset \cdots \supset X_{k}$ such that $X_{i}\in |\mathcal{L}_{i-1}|$ for $1\leq i\leq k$, where $X_{0}:=X$, $\mathcal{L}_{0}:=\mathcal{L}$ and $\mathcal{L}_{i}:=\mathcal{L}|_{X_{i}}$.
Here we note that $\dim X_{j}=n-j$.
Let $r_{p,q}:H^{p}(X_{q},\mathcal{L}_{q})\to H^{p}(X_{q+1},\mathcal{L}_{q+1})$ 
be the natural map.
\end{Definition}

\begin{Theorem}\label{Theorem1}
{\rm (\cite[Propositions 2.1 and 2.3, and Theorem 2.4]{Fukuma04-2})}
Let $X$ be a projective variety of dimension $n\geq 2$ 
and let $\mathcal{L}$ be an ample line bundle on $X$.
Assume that $h^{t}(-s\mathcal{L})=0$ for every integers $t$ and $s$ with $0\leq t\leq n-1$ and $1\leq s$, and $|\mathcal{L}|$ has an $(n-i)$-ladder for an integer $i$ with $1\leq i\leq n$.
Then the $i$th sectional geometric genus has the following properties.
\begin{itemize}
\item [\rm (1)] $g_{i}(X_{j},\mathcal{L}_{j})=g_{i}(X_{j+1},\mathcal{L}_{j+1})$ for every integer $j$ with $0\leq j\leq n-i-1$. {\rm (}Here we use the notation in Definition {\rm \ref{1.5}}.{\rm )}
\item [\rm (2)]
$g_{i}(X,\mathcal{L})\geq h^{i}(\mathcal{O}_{X})$.
\end{itemize}
\end{Theorem}

In particular, if $(X,\mathcal{L})$ is a polarized manifold with $\mbox{Bs}|\mathcal{L}|=\emptyset$, then $\mathcal{L}$ has an $(n-i)$-ladder $X\supset X_{1}\supset \cdots \supset X_{n-i}$ such that each $X_{j}$ is smooth, and from Theorem \ref{Theorem1} (1) and Remark \ref{Remark1} (2) we see that $g_{i}(X,\mathcal{L})=g_{i}(X_{n-i},\mathcal{L}_{n-i})=h^{i}(\mathcal{O}_{X_{n-i}})=h^{0}(\Omega_{X_{n-i}})$, that is,
the $i$th sectional geometric genus is the geometric genus of $i$-dimensional projective variety $X_{n-i}$.
This is a reason why we call this invariant 
the $i$th sectional geometric genus.
\par
From Theorem \ref{Theorem1} we see that the $i$th sectional geometric genus is expected to have properties similar to those of the geometric genus of $i$-dimensional projective varieties.
In particular, if $i=2$, then $g_{2}(X,\mathcal{L})$ is expected to have properties similar to those of the geometric genus of projective surfaces and we can propose several problems which can be considered as a generalization of theorems in the theory of surfaces.
In \cite{Fukuma05-2}, we investigated them. See \cite{Fukuma05-2} for further detail.
\par
For other results concerning the $i$th sectional geometric genus, for example, see \cite{Fukuma04}, \cite{Fukuma04-2}, \cite{Fukuma05-1} and \cite{Fukuma05-2}.
\end{Remark}
\noindent
\par
The following result will be used later.

\begin{Theorem}\label{1.2.1}
Let $X$ be a projective variety with $\dim X=n$ and let $\mathcal{L}$ be a nef and big line bundle on $X$. 
\\
{\rm (1)} For any integer $i$ with $0\leq i\leq n-1$, we have
$$g_{i}(X,\mathcal{L})
=\sum_{j=0}^{n-i-1}(-1)^{n-j}{n-i\choose j}\chi(-(n-i-j)\mathcal{L})
+\sum_{k=0}^{n-i}(-1)^{n-i-k}h^{n-k}(\mathcal{O}_{X}).$$
\noindent
{\rm (2)} Assume that $X$ is smooth. 
Then for any integer $i$ with $0\leq i\leq n-1$, we have
$$g_{i}(X,\mathcal{L})
=\sum_{j=0}^{n-i-1}(-1)^{j}{n-i\choose j}h^{0}(K_{X}+(n-i-j)\mathcal{L})
+\sum_{k=0}^{n-i}(-1)^{n-i-k}h^{n-k}(\mathcal{O}_{X}).$$
\end{Theorem}
\noindent{\em Proof.}
(1) By the same argument as in the proof of \cite[Theorem 2.2]{Fukuma04}, we obtain
\begin{eqnarray*}
\chi_{n-i}(X,\mathcal{L})
&=&\sum_{j=0}^{n-i}(-1)^{n-i-j}{n-i\choose j}\chi(-(n-i-j)\mathcal{L}) \\
&=&\sum_{j=0}^{n-i-1}(-1)^{n-i-j}{n-i\choose j}\chi(-(n-i-j)\mathcal{L})
+\chi(\mathcal{O}_{X}).
\end{eqnarray*}
Hence by Definition \ref{1.2}, we get the assertion. 
\\
(2) By using the Serre duality and the Kawamata-Viehweg vanishing theorem, we get the assertion from (1). $\Box$

\begin{Proposition}\label{PR3.1}
Let $(X,\mathcal{L})$ be a polarized manifold of dimension $n\geq 4$ and let $(M,\mathcal{A})$ be a reduction of $(X,\mathcal{L})$.
If $\kappa(K_{X}+(n-3)\mathcal{L})=-\infty$, then $g_{j}(X,\mathcal{L})=g_{j}(M,\mathcal{A})=0$ for every integer $j\geq 3$ unless $(M,\mathcal{A})$ is a scroll over a normal projective variety of dimension $3$.
\end{Proposition}
\noindent{\em Proof.}
Assume that $(M,\mathcal{A})$ is not a scroll over a normal projective variety of dimension $3$.
Then by Lemma \ref{LM1} we have $h^{j}(\mathcal{O}_{X})=h^{j}(\mathcal{O}_{M})=0$ for every integer $j\geq 3$.
By assumption, we get $h^{0}(K_{M}+t\mathcal{A})=h^{0}(K_{X}+t\mathcal{L})=0$ for every integer $t$ with $1\leq t\leq n-3$.
Hence by Theorem \ref{1.2.1} (2) and Remark \ref{Remark1} (5) we get $g_{j}(X,\mathcal{L})=g_{j}(M,\mathcal{A})=0$ for every integer $j\geq 3$.
This completes the proof. $\Box$

\begin{Remark}\label{TR2-3.1}
If $(M,\mathcal{A})$ is a scroll over a normal projective variety of dimension $3$,
then by \cite[Example 2.10 (8)]{Fukuma04} we have $g_{j}(X,\mathcal{L})=g_{j}(M,\mathcal{A})=0$ for every integer $j\geq 4$ and $g_{3}(X,\mathcal{L})=g_{3}(M,\mathcal{A})=h^{3}(\mathcal{O}_{M})=h^{3}(\mathcal{O}_{X})$.
\end{Remark}
\noindent
\par
As the next step, we want to generalize the notion of the $\Delta$-genus.
Several generalizations can be considered from various point of view.
Here we will give a generalization of the $\Delta$-genus from the following point of view.
For the case of $\Delta(X,\mathcal{L})$, the following result has been obtained.

\begin{Theorem}\label{Theorem2}
{\rm (}See e.g. {\rm \cite[$\S$3 in Chapter I]{FujitaBook}}.{\rm )}
Let $X$ be a projective variety of dimension $n\geq 2$ 
and let $\mathcal{L}$ be an ample line bundle on $X$.
We use the notation in Definition {\rm \ref{1.5}}.
If $|\mathcal{L}|$ has an $(n-1)$-ladder and $h^{0}(\mathcal{L}_{n-1})>0$, then 
$$\Delta(X,\mathcal{L})=\sum_{j=0}^{n-1}\dim\mbox{\rm Coker}(r_{0,j}).$$
In particular, we have $\Delta(X,\mathcal{L})\geq \Delta(X_{1},\mathcal{L}_{1})\geq \cdots \geq \Delta(X_{n-1},\mathcal{L}_{n-1})\geq 0$.
\end{Theorem}

Here we want to give the definition of the {\it $i$th $\Delta$-genus} which satisfies a generalization of Theorem \ref{Theorem2}.
Now we are going to give the definition of the $i$th $\Delta$-genus.

\begin{Definition}\label{1.3}
(\cite[Definition 2.1]{Fukuma05})
Let $(X,\mathcal{L})$ be a polarized variety of dimension $n$.
For every integer $i$ with $0\leq i\leq n$, {\it the $i$th $\Delta$-genus $\Delta_{i}(X,\mathcal{L})$ of} $(X,\mathcal{L})$ is defined by the following formula:
$$
{\Delta}_{i}(X,\mathcal{L})= 
\left\{
\begin{array}{ll}
0 & \mbox{if $i=0$,} \\
g_{i-1}(X,\mathcal{L})-{\Delta}_{i-1}(X,\mathcal{L}) \\
\ \ \ +(n-i+1)h^{i-1}(\mathcal{O}_{X})-h^{i-1}(\mathcal{L})
& \mbox{if $1\leq i\leq n$.}
\end{array}\right.
$$
\end{Definition}

\begin{Remark}\label{Remark1.3}
\begin{itemize}
\item [\rm (1)]
If $i=1$, then ${\Delta}_{1}(X,\mathcal{L})$ is equal to the $\Delta$-genus of $(X,\mathcal{L})$. 
\item [\rm (2)]
If $i=n$, then ${\Delta}_{n}(X,\mathcal{L})=h^{n}({\cal O}_{X})-h^{n}(\mathcal{L})$ (see \cite[Proposition 2.4]{Fukuma05}).
\item [\rm (3)] For every integer $i$ with $1\leq i\leq n$, by the definition of the $i$th $\Delta$-genus, we have the following equality which will be used later.
$$
\Delta_{i-1}(X,\mathcal{L})=
g_{i-1}(X,\mathcal{L})-{\Delta}_{i}(X,\mathcal{L})+(n-i+1)h^{i-1}(\mathcal{O}_{X})-h^{i-1}(\mathcal{L}).
$$
\item [\rm (4)]
Let $(X,\mathcal{L})$ be a polarized manifold of dimension $n$ and 
let $(M,\mathcal{A})$ be a reduction of $(X,\mathcal{L})$.
Then $\Delta_{i}(X,\mathcal{L})=\Delta_{i}(M,\mathcal{A})$ for every integer $i$ with $2\leq i\leq n$ (see \cite[Corollary 2.11]{Fukuma05}).
\end{itemize}
\end{Remark}
\noindent
\\
Then, for the case of the $i$th $\Delta$-genus, we can prove the following.

\begin{Theorem}\label{Theorem3}{\rm (}See {\rm \cite[Theorem 2.8 and Corollary 2.9]{Fukuma05}} and {\rm \cite[Proposition 2.1]{Fukuma04-2}}.{\rm )}
Let $X$ be a projective variety of dimension $n\geq 2$ 
and let $\mathcal{L}$ be an ample line bundle on $X$.
We use the notation in Definition {\rm \ref{1.5}}.
Assume that $h^{t}(-s\mathcal{L})=0$ for every integers $t$ and $s$ with $0\leq t\leq n-1$ and $1\leq s$.
If $|\mathcal{L}|$ has an $(n-i)$-ladder and $h^{0}(\mathcal{L}_{n-i})>0$ for an integer $i$ with $1\leq i\leq n$, then 
$$\Delta_{i}(X,\mathcal{L})=\sum_{j=0}^{n-i}\dim\mbox{\rm Coker}(r_{i-1,j}).$$
In particular, we have $\Delta_{i}(X,\mathcal{L})\geq \Delta_{i}(X_{1},\mathcal{L}_{1})\geq \cdots \geq \Delta_{i}(X_{n-i},\mathcal{L}_{n-i})\geq 0$.
\end{Theorem}

The definition of the $i$th $\Delta$-genus is so complicated that
a lot of things about the $i$th $\Delta$-genus are unknown.
So it is important to investigate the following problems in order to 
understand the meaning and properties of the $i$th $\Delta$-genus.

\begin{Problem}\label{Problem}
\begin{itemize}
\item [\rm (i)] Does the $i$th $\Delta$-genus have properties similar to those of the $\Delta$-genus $?$ For example, there are the following two questions.
\begin{itemize}
\item [\rm (i.1)] Does $\Delta_{i}(X,\mathcal{L})\geq 0$ hold $?$
\item [\rm (i.2)] Can we get the $i$th $\Delta$-genus version of the Fujita theory on $\Delta$-genus $?$
\end{itemize}
\item [\rm (ii)] Are there any relationship 
between $g_{i}(X,\mathcal{L})$ and $\Delta_{i}(X,\mathcal{L})$ $?$
\item [\rm (iii)] Are there any relationship 
between $\Delta_{i}(X,\mathcal{L})$ and $\Delta_{i+1}(X,\mathcal{L})$ $?$
\item [\rm (iv)] Classify $(X,\mathcal{L})$ by the value of the $i$th $\Delta$-genus.
\item [\rm (v)] What is the geometric meaning of the $i$th $\Delta$-genus $?$
\end{itemize}
\end{Problem}

\begin{Remark}\label{TR1}
If $X$ is smooth and $\mathcal{L}$ is ample, then the following facts on Problem \ref{Problem} are known.
\begin{itemize}
\item [\rm (1)] 
First we consider Problem \ref{Problem} (i.1).
If $i=1$, then $\Delta_{1}(X,\mathcal{L})\geq 0$ (see \cite[(4.2) Theorem]{FujitaBook}).
Moreover if $\mathcal{L}$ is base point free, then $\Delta_{i}(X,\mathcal{L})\geq 0$ holds for every integer $i$ with $0\leq i\leq n$.
But unfortunately, there exists an example $(X,\mathcal{L})$ such that $\Delta_{i}(X,\mathcal{L})<0$ (see \cite[Section 4]{Fukuma05}).
\item [\rm (2)]
Next we consider Problem \ref{Problem} (ii).
If $i=1$ and $\mathcal{L}$ is merely ample, then it is known that 
$g_{1}(X,\mathcal{L})=0$ if and only if $\Delta_{1}(X,\mathcal{L})=0$ (see \cite[(12.1) Theorem]{FujitaBook}).
Next we consider the case of $i\geq 2$.
Then under the assumption that $\mbox{Bs}|\mathcal{L}|=\emptyset$ we see that $g_{i}(X,\mathcal{L})=0$ if and only if $\Delta_{i}(X,\mathcal{L})=0$ (see \cite[Theorem 3.13]{Fukuma05}).
\item [\rm (3)]
Next we consider Problem \ref{Problem} (iii) under the assumption that $\mathcal{L}$ is base point free.
Then, for example, we get the following:
If $\Delta_{i}(X,\mathcal{L})\leq i-1$, then $\Delta_{i+1}(X,\mathcal{L})=0$ (see \cite[Proposition 3.9]{Fukuma05}).
In particular, if $\Delta_{i}(X,\mathcal{L})=0$, then $\Delta_{i+1}(X,\mathcal{L})=0$.
Maybe there will be several relationship between $\Delta_{i}(X,\mathcal{L})$ and $\Delta_{i+1}(X,\mathcal{L})$ other than this.
\item [\rm (4)]
For Problem \ref{Problem} (iv), there exists the classification of $(X,\mathcal{L})$ by the value of $\Delta_{2}(X,\mathcal{L})$ as follows:
\begin{itemize}
\item [\rm (4.1)] The classification of polarized manifolds $(X,\mathcal{L})$ such that $\mbox{\rm Bs}|\mathcal{L}|=\emptyset$ and $\Delta_{2}(X,\mathcal{L})=0$ (see \cite[Theorem 3.13 and Remark 3.13.1]{Fukuma05}).
\item [\rm (4.2)] The classification of polarized manifolds $(X,\mathcal{L})$ such that $\mathcal{L}$ is very ample and $\Delta_{2}(X,\mathcal{L})=1$ (see \cite[Theorem 3.17]{Fukuma05} and \cite[Remark 2]{Fukuma08-4}).
\end{itemize}

\item [\rm (5)]
At present, we do not know much about any answer to Problem \ref{Problem} (v).
This problem seems to be the most difficult problem among the above problems even in the case where $\mathcal{L}$ is base point free or very ample.
\end{itemize}
\end{Remark}

In this paper, we consider Problem \ref{SGG-PB1} (i.2) and Problem \ref{Problem} (i.2).
In \cite[(1.9) Theorem]{Fujita80}, Fujita proved that $(X,\mathcal{L})$ is a Del Pezzo manifold (namely $K_{X}=-(n-1)\mathcal{L}$) if and only if $g(X,\mathcal{L})=1$ and $\Delta(X,\mathcal{L})=1$, that is, $g_{1}(X,\mathcal{L})=1$ and $\Delta_{1}(X,\mathcal{L})=1$.
So in this paper, we consider an analogous characterization of $(X,\mathcal{L})$ with $K_{X}=-(n-i)\mathcal{L}$ by using $g_{i}(X,\mathcal{L})$ and $\Delta_{i}(X,\mathcal{L})$.

\section{Main Theorems}\label{S4}

\subsection{A conjecture}

First we provide the following conjecture which is the main theme of this paper.
\begin{s-Conjecture}\label{CJ}
Let $(X,\mathcal{L})$ be a polarized manifold of dimension $n\geq 3$.
Then, for every integer $i$ with $2\leq i\leq n-1$, the following are equivalent one another.
\par
{\rm $C(i,1)$:} $K_{X}=-(n-i)\mathcal{L}$.
\par
{\rm $C(i,2)$:} $\Delta_{i}(X,\mathcal{L})=1$ and $2g_{1}(X,\mathcal{L})-2=(i-1)\mathcal{L}^{n}$.
\par
{\rm $C(i,3)$:} $\Delta_{i}(X,\mathcal{L})>0$ and $2g_{1}(X,\mathcal{L})-2=(i-1)\mathcal{L}^{n}$.
\par
{\rm $C(i,4)$:} $g_{i}(X,\mathcal{L})=1$ and $2g_{1}(X,\mathcal{L})-2=(i-1)\mathcal{L}^{n}$.
\par
{\rm $C(i,5)$:} $g_{i}(X,\mathcal{L})>0$ and $2g_{1}(X,\mathcal{L})-2=(i-1)\mathcal{L}^{n}$.
\end{s-Conjecture}

\begin{s-Remark}\label{MT-Rem}
If $i=1$, then $C(1,1)$ and $C(1,2)$ in Conjecture \ref{CJ} are equivalent each other for {\it any} ample line bundle $\mathcal{L}$.
Of course $C(1,1)$ implies $C(1,3)$ (resp. $C(1,4)$, $C(1,5)$).
But $C(1,3)$ (resp. $C(1,4)$, $C(1,5)$) does not imply $C(1,1)$ because $(X,\mathcal{L})$ is possibly a scroll over an elliptic curve.
\end{s-Remark}

\begin{s-Remark}\label{MT-Rem2}
As a generalization of the case where $i=1$, it is natural to consider that $C(i,1)$ is equivalent to the following condition:
\par
{\rm $C(i,6)$:} $\Delta_{i}(X,\mathcal{L})=1$ and $g_{i}(X,\mathcal{L})=1$.
\\
We can easily see that $C(i,1)$ implies $C(i,6)$.
But from the following examples (Examples \ref{EX1} and \ref{EX2}) we see that its converse is not true in general.
\end{s-Remark}

\begin{s-Example}\label{EX1}
Let $n$ be a natural number with $n\geq 3$, and 
let $Y$ be a smooth projective variety of dimension $m$ 
with $1\leq m\leq n-2$.
Let $\mathcal{H}$ be an ample line bundle on $Y$ such that 
$K_{Y}\neq -(n-m-1)\mathcal{H}$
and $h^{0}(K_{Y}+(n-m-1)\mathcal{H})=1$.
(There exists a polarized manifold $(Y,\mathcal{H})$ like this.
For example, let $Y$ be a principally polarized Abelian variety with $\dim Y=m=n-2$ and let $\mathcal{H}$ be an ample line bundle on $Y$ such that $\mathcal{H}^{m}=m!$.
Then $K_{Y}+(n-m-1)\mathcal{H}=\mathcal{H}$ and 
$h^{0}(K_{Y}+(n-m-1)\mathcal{H})=h^{0}(\mathcal{H})=1$.)
\par
Next we take a Del Pezzo manifold $(F,\mathcal{A})$ of dimension $n-m$.
Then we note that $K_{F}=-(n-m-1)\mathcal{A}$.
\par
Here we set $X:=Y\times F$ and $L:=p_{1}^{*}(\mathcal{H})\otimes p_{2}^{*}(\mathcal{A})$,
where $p_{i}$ denotes the $i$th projection map.
Then $K_{X}+(n-m-1)\mathcal{L}=p_{1}^{*}(K_{Y}+(n-m-1)\mathcal{H})$.
By \cite[Lemma 1.6]{Fukuma05} we also get $h^{j}(\mathcal{O}_{X})=0$ and $h^{j}(\mathcal{L})=0$ for every integer $j$ with $j\geq m+1$.
Hence $\Delta_{n}(X,\mathcal{L})=0$ by Remark \ref{Remark1.3} (2), and by Theorem \ref{1.2.1} (2)
we see that 
$g_{j}(X,\mathcal{L})=0$ for every integer $j$ with $j\geq m+2$ and 
$g_{m+1}(X,\mathcal{L})=h^{0}(K_{X}+(n-m-1)\mathcal{L})=1$.
Moreover by Remark \ref{Remark1.3} (3) we see that
$\Delta_{j}(X,\mathcal{L})=0$ for every integer $j\geq m+2$ and
\begin{eqnarray*}
&&\Delta_{m+1}(X,\mathcal{L})\\
&&=g_{m+1}(X,\mathcal{L})-\Delta_{m+2}(X,\mathcal{L})+(n-m-1)h^{m+1}(\mathcal{O}_{X})-h^{m+1}(\mathcal{L}) \\
&&=1.
\end{eqnarray*}
Therefore $g_{m+1}(X,\mathcal{L})=\Delta_{m+1}(X,\mathcal{L})=1$.
But $K_{X}\neq-(n-m-1)\mathcal{L}$ and this $(X,\mathcal{L})$ is an example.
\end{s-Example}

\begin{s-Example}\label{EX2}
Let $k$ be a natural number with $k\geq 2$ 
and set $n:=2k+1$ and $i:=(n-1)/2$.
Here we consider $(M,\mathcal{A})=(\mathbb{P}^{n}, \mathcal{O}_{\mathbb{P}^{n}}(2))$.
Then $K_{M}=-(k+1)\mathcal{A}=-(n-i)\mathcal{A}$.
Moreover we see that 
$g_{i}(M,\mathcal{A})=1$ and $\Delta_{i}(M,\mathcal{A})=1$ (see (I) in the proof of Theorem \ref{MT4} below). 
Let $\pi:X\to \mathbb{P}^{n}$ be the blowing up at a general point on $\mathbb{P}^{n}$ and let $\mathcal{L}:=\pi^{*}(\mathcal{A})-E$, where $E$ is the exceptional divisor.
Then by \cite[Theorem 0.1]{Ballico99}, we see that $(X,\mathcal{L})$ is a polarized manifold with $K_{X}+(n-i)\mathcal{L}=(i-1)E$.
On the other hand, we note that $(M,\mathcal{A})$ is a reduction of $(X,\mathcal{L})$ and $2\leq i<n-1$.
Hence by Remarks \ref{Remark1} (5) and \ref{Remark1.3} (4) we get $g_{i}(X,\mathcal{L})=g_{i}(M,\mathcal{A})=1$ and $\Delta_{i}(X,\mathcal{L})=\Delta_{i}(M,\mathcal{A})=1$.
\end{s-Example}

\subsection{The case where \textrm{\boldmath $\mbox{max}\{ 2, \dim\mbox{Bs}|\mathcal{L}|+2 \}\leq i \leq n-1$}}

First we consider the case where $\mbox{\rm max}\{ 2, \dim\mbox{Bs}|\mathcal{L}|+2 \}\leq i\leq n-1$.

\begin{s-Theorem}\label{MT4}
Let $(X,\mathcal{L})$ be a polarized manifold of dimension $n\geq 3$, and let $m=\dim \mbox{\rm Bs}|\mathcal{L}|$. {\rm (}If $\mbox{\rm Bs}|\mathcal{L}|=\emptyset$, then we set $m=-1$.{\rm)}
Assume that $i$ is an integer with $\mbox{\rm max}\{ 2, m+2 \}\leq i \leq n-1$.
Then Conjecture {\rm \ref{CJ}} is true.
\end{s-Theorem}
\noindent{\em Proof.}
By assumption and \cite[Proposition 1.12 (2)]{Fukuma04-2}, we see that the following hold:
\begin{itemize}
\item [$(A_{i})$] $\mathcal{L}$ has an $(n-i)$-ladder $X_{n-i}\subset \cdots \subset X_{1}\subset X$.
\item [$(B_{i})$] $h^{0}(\mathcal{L}_{n-i})>0$.
\item [(C)] $h^{j}(\mathcal{L}^{\otimes -t})=0$ for any $j$ and $t$ with $0\leq j\leq n-1$ and $t>0$.
\item [$(D_{i})$] $X_{j}$ is normal for any $j$ with $0\leq j\leq n-i$.
\item [$(E_{i})$] $X_{j}$ is Cohen-Macaulay for any $j$ with $0\leq j\leq n-i$.
\end{itemize}
\noindent
\\
(I) Assume that $C(i,1)$ holds.
Then by Remark \ref{Remark1} (4) we see that 
\begin{eqnarray*}
2g_{1}(X,\mathcal{L})-2
&=&(K_{X}+(n-i)\mathcal{L}+(i-1)\mathcal{L})\mathcal{L}^{n-1}\\
&=&(i-1)\mathcal{L}^{n}.
\end{eqnarray*}
Here we note that $h^{j}(\mathcal{O}_{X})=0$ and $h^{j}(\mathcal{L})=0$ for every integer $j$ with $j\geq 2$.
Moreover $h^{0}(K_{X}+(n-i)\mathcal{L})=1$ and $h^{0}(K_{X}+k\mathcal{L})=0$ for every integer $k$ with $1\leq k\leq n-i-1$.
Hence we see that $g_{i}(X,\mathcal{L})=1$ and $g_{k}(X,\mathcal{L})=0$ for every integer $k$ with $k\geq i+1$ by Theorem \ref{1.2.1} (2) (this means that $C(i,1)$ implies $C(i,4)$), and
by Remark \ref{Remark1.3} (2) and (3), we have
$\Delta_{k}(X,\mathcal{L})=0$ for every integer $k$ with $k\geq i+1$ and
\begin{eqnarray*}
\Delta_{i}(X,\mathcal{L})
&=&g_{i}(X,\mathcal{L})-\Delta_{i+1}(X,\mathcal{L})+(n-i)h^{i}(\mathcal{O}_{X})-h^{i}(\mathcal{L})\\
&=&1.
\end{eqnarray*}
Therefore we see that $C(i,1)$ implies $C(i,2)$ and $C(i,4)$ above.
\\
\\
(II) It is trivial that $C(i,2)$ implies $C(i,3)$, and $C(i,4)$ implies $C(i,5)$.
\\
\\
(III) Assume that $C(i,3)$ holds.
Then we will prove that $C(i,5)$ holds.
In order to prove $C(i,5)$, it suffices to show that $g_{i}(X,\mathcal{L})>0$.
Here we note that $g_{i}(X,\mathcal{L})\geq 0$ by \cite[Theorem 2.4]{Fukuma04-2}.
Assume that $g_{i}(X,\mathcal{L})=0$.
Here we prove the following.
\begin{s-Claim}\label{T-CL1}
If $g_{i}(X,\mathcal{L})=0$, then $\Delta_{i}(X,\mathcal{L})=0$.
\end{s-Claim}
\noindent
{\em Proof.}
Assume that $g_{i}(X,\mathcal{L})=0$.
Then $0=g_{i}(X,L)=g_{i}(X_{n-i},\mathcal{L}_{n-i})=h^{i}(\mathcal{O}_{X_{n-i}})$ by Theorem \ref{Theorem1} and Remark \ref{Remark1} (2).
Therefore $h^{i}(\mathcal{O}_{X})=h^{i}(\mathcal{O}_{{X}_{1}})=\cdots =h^{i}(\mathcal{O}_{{X}_{n-i-1}})\leq h^{i}(\mathcal{O}_{{X}_{n-i}})=0$ by \cite[Proposition 2.1 (b)]{Fukuma04-2}.
Hence $H^{i-1}(\mathcal{L}_{j})\to H^{i-1}(\mathcal{L}_{j+1})$ is surjective for $0\leq j \leq n-i$.
Namely $\dim\mbox{Coker}(r_{i-1,j})=0$ for $0\leq j\leq n-i$.
On the other hand by Theorem \ref{Theorem3}, we have
$$
{\Delta}_{i}(X,\mathcal{L})=\sum_{k=0}^{n-i}\dim\mbox{Coker}(r_{i-1,k}).
$$

Therefore we get 
$$
{\Delta}_{i}(X,\mathcal{L})=\sum_{k=0}^{n-i}\dim\mbox{Coker}(r_{i-1,k})=0.
$$
This completes the proof of Claim \ref{T-CL1}. $\Box$
\\
\par
But this contradicts the assumption.
Therefore $g_{i}(X,\mathcal{L})>0$ and we see that $C(i,5)$ holds.
\\
\\
(IV) Assume that $C(i,5)$ holds.
Then
\begin{eqnarray*}
1+\frac{1}{2}(i-1)\mathcal{L}^{n}
&=&g_{1}(X,\mathcal{L}) \\
&=&1+\frac{1}{2}(K_{X}+(n-1)\mathcal{L})\mathcal{L}^{n-1}\\
&=&1+\frac{1}{2}(K_{X_{n-i}}+(i-1)\mathcal{L}_{n-i}){\mathcal{L}_{n-i}}^{i-1}\\
&=&1+\frac{1}{2}(i-1)\mathcal{L}^{n}+\frac{1}{2}K_{X_{n-i}}\mathcal{L}_{n-i}^{i-1}.
\end{eqnarray*}
Hence $K_{X_{n-i}}\mathcal{L}_{n-i}^{i-1}=0$.
On the other hand, we get $g_{i}(X,\mathcal{L})=h^{i}(\mathcal{O}_{X_{n-i}})$ by $(A_{i})$ and (C) (see also \cite[Propositions 2.1 and 2.3]{Fukuma04-2}.
Furthermore by $(D_{i})$, $(E_{i})$ and the Serre duality, we obtain $h^{0}(K_{X_{n-i}})=h^{i}(\mathcal{O}_{X_{n-i}})$.
Hence we have $0<g_{i}(X,\mathcal{L})=h^{i}(\mathcal{O}_{X_{n-i}})=h^{0}(K_{X_{n-i}})$.
Hence we see that $K_{X_{n-i}}=\mathcal{O}_{X_{n-i}}$.
\\
\par
Next we prove the following claim.
\begin{s-Claim}\label{CL2}
A natural map $\mbox{\rm Pic}(X_{j})\to \mbox{\rm Pic}(X_{j+1})$ is injective for any $j$ with $0\leq j\leq n-i-1$.
\end{s-Claim}
\noindent
{\em Proof.}
From the following exact sequence
$$0\to \mathbb{Z}\to \mathcal{O}_{X_{j}}\to \mathcal{O}_{X_{j}}^{*}\to 0,$$
we get the following commutative diagram.

\setlength{\unitlength}{1mm}
\begin{picture}(150,40)(-10,0)
\put(0,0){\makebox(10,10)[r]{$H^{1}(X_{j+1},\mathbb{Z})$}}
\put(30,0){\makebox(10,10)[r]{$H^{1}(\mathcal{O}_{X_{j+1}})$}}
\put(60,0){\makebox(10,10)[r]{$H^{1}(\mathcal{O}_{X_{j+1}}^{*})$}}
\put(90,0){\makebox(10,10)[r]{$H^{2}(X_{j+1},\mathbb{Z})$}}
\put(0,25){\makebox(10,10)[r]{$H^{1}(X_{j},\mathbb{Z})$}}
\put(30,25){\makebox(10,10)[r]{$H^{1}(\mathcal{O}_{X_{j}})$}}
\put(60,25){\makebox(10,10)[r]{$H^{1}(\mathcal{O}_{X_{j}}^{*})$}}
\put(90,25){\makebox(10,10)[r]{$H^{2}(X_{j},\mathbb{Z})$}}

\put(12,5){\vector(1,0){10}}
\put(42,5){\vector(1,0){10}}
\put(71,5){\vector(1,0){9}}
\put(13,30){\vector(1,0){10}}
\put(45,30){\vector(1,0){10}}
\put(72,30){\vector(1,0){10}}
\put(0,25){\vector(0,-1){15}}
\put(30,25){\vector(0,-1){15}}
\put(60,25){\vector(0,-1){15}}
\put(90,25){\vector(0,-1){15}}
\put(-5,18){$\varphi_{1}$}
\put(25,18){$\varphi_{2}$}
\put(55,18){$\varphi_{3}$}
\put(85,18){$\varphi_{4}$}

\end{picture}

So in order to prove that $\varphi_{3}$ is injective, it suffices to show the following for every integer $j$ with $0\leq j\leq n-i-1$ because $\mbox{Pic}(X_{j})\cong H^{1}(\mathcal{O}_{X_{j}}^{*})$ and $\mbox{Pic}(X_{j+1})\cong H^{1}(\mathcal{O}_{X_{j+1}}^{*})$.
\begin{itemize}
\item [(a)] $h^{1}(\mathcal{O}_{X_{j}}(-X_{j+1}))=0$.
\item [(b)] $H^{1}(X_{j},\mathbb{Z})\cong H^{1}(X_{j+1}, \mathbb{Z})$.
\item [(c)] The map $H^{2}(X_{j},\mathbb{Z})\to H^{2}(X_{j+1}, \mathbb{Z})$ is injective.
\end{itemize}

By (C) we can prove $h^{t}(\mathcal{L}_{j}^{\otimes -s})=0$ for every $j$, $t$ and $s$ with $0\leq j\leq n-i-1$, $0\leq t\leq n-j-1$ and $1\leq s$.
Therefore we get (a) since $\dim X_{n-i-1}=i+1\geq 3$.
\par
Next we consider (b) and (c).
In this case we need to take an $(n-i)$-ladder carefully.
Namely, we take general members $X_{1}\in |\mathcal{L}|, X_{2}\in |\mathcal{L}|_{X_{1}}, \dots , X_{n-i}\in |\mathcal{L}|_{X_{n-i-1}}$.
Then $X\supset X_{1}\supset \dots \supset X_{n-i}$ is an $(n-i)$-ladder such that $X_{j}-X_{j+1}$ is smooth for every $j$.
Hence $X_{j}-X_{j+1}$ is locally complete intersection.
Here we use \cite[Corollary 2.3.3]{BeSoBook}.
Since $X_{j+1}$ is an ample line bundle on $X_{j}$, we see that if $2\leq i=\dim X_{n-i-1}-1$, then $H^{t}(X_{j},\mathbb{Z})\to H^{t}(X_{j+1}, \mathbb{Z})$ is an isomorphism (resp. injective) for $t=1$ (resp. $t=2$) and every $j$ with $0\leq j\leq n-i-1$.
\par
Therefore we get the assertion of Claim \ref{CL2}. $\Box$
\\
\par
By this claim we have $K_{X}+(n-i)\mathcal{L}=\mathcal{O}_{X}$.
So we get $C(i,1)$.
This completes the proof of Theorem \ref{MT4}. $\Box$

\subsection{The case of \textrm{\boldmath $i=2$}}

Next we consider the case where $i=2$ and $\mathcal{L}$ is ample in general.
Then we can prove the following:

\begin{s-Theorem}\label{MT2}
Let $(X,\mathcal{L})$ be a polarized manifold of dimension $n\geq 3$.
Then Conjecture {\rm \ref{CJ}} for $i=2$ is true.
\end{s-Theorem}
\noindent{\em Proof.}
(I) First we assume that $C(2,1)$. Then by the same argument as in the proof of Theorem \ref{MT4} we see that $C(2,2)$ and $C(2,4)$ hold.
\\
\\
(II) It is trivial that $C(2,2)$ (resp. $C(2,4)$) implies $C(2,3)$ (resp. $C(2,5)$).
\\
\\
(III) Assume that $C(2,3)$.
Then we will prove that $C(2,5)$ holds.
In order to prove $C(2,5)$, it suffices to show that $g_{2}(X,\mathcal{L})>0$.
By the assumption that $2g_{1}(X,\mathcal{L})-2=\mathcal{L}^{n}$, we get $(K_{X}+(n-2)\mathcal{L})\mathcal{L}^{n-1}=0$.
Hence by \cite[Lemma 2.5.9]{BeSoBook} we have $\kappa(K_{X}+(n-2)\mathcal{L})\leq 0$.
\\
(III.1) If $\kappa(K_{X}+(n-2)\mathcal{L})=-\infty$, then $(X,\mathcal{L})$ is one of the types (1) to (7.4) in Theorem \ref{TH1} by Remark \ref{RM2} (1).
By \cite[Example 2.10]{Fukuma04} and \cite[Example 2.12]{Fukuma05} we may assume that $(X,\mathcal{L})$ is a scroll over a smooth surface $S$ because we assume that $\Delta_{2}(X,\mathcal{L})>0$.
In this case, 
by \cite[(3.2.1)]{BeSoWi} and \cite[(11.8.6) in the proof of (11.8) Theorem]{FujitaBook}, there exists an ample vector bundle $\mathcal{E}$ of rank $n-1$ on $X$ such that $X=\mathbb{P}_{S}(\mathcal{E})$, $\mathcal{L}=H(\mathcal{E})$.
Let $\pi: X\to S$ be its morphism.
Here we calculate $(K_{X}+(n-2)\mathcal{L})\mathcal{L}^{n-1}$.
\begin{eqnarray*}
(K_{X}+(n-2)\mathcal{L})\mathcal{L}^{n-1}
&=&(-H(\mathcal{E})+\pi^{*}(K_{S}+c_{1}(\mathcal{E}))H(\mathcal{E})^{n-1}\\
&=&K_{S}c_{1}(\mathcal{E})+c_{2}(\mathcal{E}).
\end{eqnarray*}
If $h^{2}(\mathcal{O}_{S})=0$, 
then $h^{2}(\mathcal{O}_{X})=h^{2}(\mathcal{O}_{S})=0$ and 
$\Delta_{2}(X,\mathcal{L})=(n-1)h^{2}(\mathcal{O}_{X})-h^{2}(\mathcal{L})=-h^{2}(\mathcal{L})\leq 0$ and this contradicts the assumption.
Hence $h^{2}(\mathcal{O}_{S})\geq 1$ and by the Serre duality we have $h^{0}(K_{S})\geq 1$.
Since $\mathcal{E}$ is ample, we see that $K_{S}c_{1}(\mathcal{E})\geq 0$ and $c_{2}(\mathcal{E})>0$.
Therefore $(K_{X}+(n-2)\mathcal{L})\mathcal{L}^{n-1}>0$ and this contradicts the assumption.
Therefore there does not exist any $(X,\mathcal{L})$ with $\kappa(K_{X}+(n-2)\mathcal{L})=-\infty$, $\Delta_{2}(X,\mathcal{L})>0$ and $2g_{1}(X,\mathcal{L})-2=\mathcal{L}^{n}$.
\\
(III.2) Assume that $\kappa(K_{X}+(n-2)\mathcal{L})=0$.
Let $(M,\mathcal{A})$ be a reduction of $(X,\mathcal{L})$.
Then by Theorem \ref{TH1}, $(M,\mathcal{A})$ is a Mukai manifold and by \cite[Proposition 2.6 and Example 2.10 (7)]{Fukuma04} we see that $g_{2}(X,\mathcal{L})=g_{2}(M,\mathcal{A})>0$.
Hence by (III.1) and (III.2) we get $C(2,5)$.
\\
\\
(IV) Assume that $C(2,5)$.
Then by the same argument as (III) above, we see that 
$\kappa(K_{X}+(n-2)\mathcal{L})\leq 0$.
If $\kappa(K_{X}+(n-2)\mathcal{L})=-\infty$, then $(X,\mathcal{L})$ is a scroll over a smooth surface $S$ since $g_{2}(X,\mathcal{L})>0$.
On the other hand $h^{2}(\mathcal{O}_{X})>0$ because $g_{2}(X,\mathcal{L})=h^{2}(\mathcal{O}_{X})$.
Hence by the same argument as (III.1) above, we see that $(K_{X}+(n-2)\mathcal{L})\mathcal{L}^{n-1}>0$ and this is impossible.
Therefore $\kappa(K_{X}+(n-2)\mathcal{L})=0$ and $K_{M}+(n-2)\mathcal{A}=\mathcal{O}_{X}$, where $(M,\mathcal{A})$ is a reduction of $(X,\mathcal{L})$.
Here we prove that $(X,\mathcal{L})\cong (M,\mathcal{A})$.
So we assume that $(X,\mathcal{L})\not\cong (M,\mathcal{A})$.
Then $(K_{X}+(n-2)\mathcal{L})\mathcal{L}^{n-1}>(K_{M}+(n-2)\mathcal{A})\mathcal{A}^{n-1}$ holds.
But since $(K_{X}+(n-2)\mathcal{L})\mathcal{L}^{n-1}=0$ and $(K_{M}+(n-2)\mathcal{A})\mathcal{A}^{n-1}=0$,
this is impossible.
Hence $(X,\mathcal{L})\cong (M,\mathcal{A})$ and we get $C(2,1)$.
\\
\par
This completes the proof of Theorem \ref{MT2}. $\Box$

\begin{s-Corollary}\label{CR1}
Let $(X,\mathcal{L})$ be a polarized manifold of dimension $n\geq 3$.
Assume that $\dim \mbox{\rm Bs}|\mathcal{L}|\leq 1$. 
Then Conjecture {\rm \ref{CJ}} is true.
\end{s-Corollary}
\noindent
{\em Proof.}
Since $\dim \mbox{Bs}|\mathcal{L}|\leq 1$, we see that Conjecture \ref{CJ} is true for $i\geq 3$ by Theorem \ref{MT4}.
On the other hand, if $i=2$, then Conjecture \ref{CJ} is also true by Theorem \ref{MT2}.
Therefore we get the assertion. $\Box$

\subsection{The case where \textrm{\boldmath $i=3$} and \textrm{\boldmath $n\geq 5$}}

Next we consider the case where $i=3$ and $n\geq 5$.

\begin{s-Theorem}\label{MT3}
Let $(X,\mathcal{L})$ be a polarized manifold of dimension $n\geq 5$.
Then Conjecture {\rm \ref{CJ}} for $i=3$ is true.
\end{s-Theorem}
{\em Proof.}
(I) 
By the same argument as in the proof of Theorem \ref{MT4}, we see that $C(3,1)$ implies $C(3,2)$ and $C(3,4)$.
\\
\\
(II) It is trivial that $C(3,2)$ (resp. $C(3,4)$) implies $C(3,3)$ (resp. $C(3,5)$).
\\
\\
(III) Assume that $C(3,3)$.
Then we will prove that $C(3,1)$ holds.
\par
Since $2g_{1}(X,\mathcal{L})-2=2\mathcal{L}^{n}$, we have $(K_{X}+(n-3)\mathcal{L})\mathcal{L}^{n-1}=0$. Hence by an argument similar to (III) in the proof of Theorem \ref{MT2}, we see that $\kappa(K_{X}+(n-3)\mathcal{L})\leq 0$.
\\
(III.1) Assume that $\kappa(K_{X}+(n-3)\mathcal{L})=0$.
Then since $\kappa(K_{X}+(n-3)\mathcal{L})=0$ and $(K_{X}+(n-3)\mathcal{L})\mathcal{L}^{n-1}=0$, there exists a positive integer $t$ such that $t(K_{X}+(n-3)\mathcal{L})=\mathcal{O}_{X}$.
But by \cite[Lemma 3.3.2]{BeSoBook} we have $K_{X}+(n-3)\mathcal{L}=\mathcal{O}_{X}$.
\\
(III.2) Assume that $\kappa(K_{X}+(n-3)\mathcal{L})=-\infty$.
\begin{s-Lemma}\label{LM3}
There does not exist any $(X,\mathcal{L})$ with $\kappa(K_{X}+(n-3)\mathcal{L})=-\infty$, $\Delta_{3}(X,\mathcal{L})>0$ and $2g_{1}(X,\mathcal{L})-2=2\mathcal{L}^{n}$.
\end{s-Lemma}
\noindent{\em Proof.}
If $\kappa(K_{X}+(n-3)\mathcal{L})=-\infty$, then by Proposition \ref{PP2} $(X,\mathcal{L})$ satisfies either (1), (2.1) or (2.2) in Proposition \ref{PP2}.
\\
(i) If $(X,\mathcal{L})$ satisfies (1) in Proposition \ref{PP2}, then by using \cite[Example 2.12]{Fukuma05} we see that $\Delta_{3}(X,\mathcal{L})=0$
unless $(M,\mathcal{A})$ is a scroll over a normal projective variety of dimension $3$.
\par
We consider the case where a reduction $(M,\mathcal{A})$ of $(X,\mathcal{L})$ is a scroll over a normal projective variety $Y$ with $\dim Y=m\geq 2$.
Then by \cite[(3.2.1) Theorem]{BeSoWi} and \cite[Proposition 2.5]{BaWi96}, we get the following:

\begin{s-Proposition}\label{PP4}
Let $(X,\mathcal{L})$ be a scroll over a $3$-dimensional normal projective variety $Y$.
If $\dim X\geq 5$, then $Y$ is smooth and $(X,\mathcal{L})$ is a classical scroll over $Y$.
\end{s-Proposition}

So we see that
$(M,\mathcal{A})$ is a classical scroll, that is, $Y$ is smooth and $(M,\mathcal{A})=(\mathbb{P}_{Y}(\mathcal{E}),H(\mathcal{E}))$, 
where $\mathcal{E}$ is an ample vector bundle on $Y$.
Then in general the following claim holds.

\begin{s-Claim}\label{CL0}
Let $(X,\mathcal{L})$ be a polarized manifold of dimension $n$ and let $(M,\mathcal{A})$ be a reduction of $(X,\mathcal{L})$.
Assume that there exists a smooth projective variety $Y$ of dimension $m\geq 2$ and an ample vector bundle $\mathcal{E}$ on $Y$ of rank $n-m+1$ such that $M=\mathbb{P}_{Y}(\mathcal{E})$ and $\mathcal{A}=H(\mathcal{E})$.
\\
{\rm (1)} If $g_{m}(X,\mathcal{L})>0$, then $2g_{1}(X,\mathcal{L})-2>(m-1)\mathcal{L}^{n}$.
\\
{\rm (2)} If $\Delta_{m}(X,\mathcal{L})>0$, then $2g_{1}(X,\mathcal{L})-2>(m-1)\mathcal{L}^{n}$.
\end{s-Claim}
{\em Proof.}
Assume that $h^{m}(\mathcal{O}_{Y})\geq 1$.
Then $\kappa(Y)\geq 0$.
On the other hand,
\begin{eqnarray*}
&&2g_{1}(X,\mathcal{L})-2-(m-1)\mathcal{L}^{n}\\
&&\geq 2g_{1}(M,\mathcal{A})-2-(m-1)\mathcal{A}^{n}\\
&&=(K_{M}+(n-m)\mathcal{A})\mathcal{A}^{n-1}\\
&&=(-H(\mathcal{E})+f^{*}(K_{Y}+\det \mathcal{E}))H(\mathcal{E})^{n-1}\\
&&=-s_{m}(\mathcal{E})+(K_{Y}+\det\mathcal{E})s_{m-1}(\mathcal{E})\\
&&=s_{m-1}(\mathcal{E})s_{1}(\mathcal{E})-s_{m}(\mathcal{E})+K_{Y}s_{m-1}(\mathcal{E}).
\end{eqnarray*}
(Here $f:M\to Y$ denotes the projection.)
Since $\mathcal{E}$ is ample, $m\geq 2$ and $\kappa(Y)\geq 0$, we have $s_{m-1}(\mathcal{E})s_{1}(\mathcal{E})-s_{m}(\mathcal{E})>0$ and $K_{Y}s_{m-1}(\mathcal{E})\geq 0$ by \cite[Example 12.1.7 and Lemma 14.5.1]{Fulton}.
Therefore we see that $2g_{1}(X,\mathcal{L})-2>(m-1)\mathcal{L}^{n}$. 
\\
(1) If $g_{m}(X,\mathcal{L})>0$, then $h^{m}(\mathcal{O}_{Y})>0$ because $g_{m}(X,\mathcal{L})=g_{m}(M,\mathcal{A})=h^{m}(\mathcal{O}_{M})=h^{m}(\mathcal{O}_{Y})$
by \cite[Example 2.10 (8)]{Fukuma04}.
Hence by the above argument we get the assertion (1).
\\
(2) Assume that $\Delta_{m}(X,\mathcal{L})>0$.
Here we note that $g_{j}(M,\mathcal{A})=0$, $h^{j}(\mathcal{O}_{M})=0$ and $h^{j}(\mathcal{A})=0$ for every $j\geq m+1$ and $g_{m}(M,\mathcal{A})=h^{m}(\mathcal{O}_{M})$ by \cite[Example 2.10 (8)]{Fukuma04} and \cite[Lemma 1.6]{Fukuma05}.
Therefore we get $\Delta_{j}(M,\mathcal{A})=0$ for every integer $j$ with $j\geq m+1$, and by using Remark \ref{Remark1.3} (3) we have $\Delta_{m}(M,\mathcal{A})=(n-m+1)h^{m}(\mathcal{O}_{M})-h^{m}(\mathcal{A})$.
If $h^{m}(\mathcal{O}_{M})=0$, then $\Delta_{m}(M,\mathcal{A})=-h^{m}(\mathcal{A})\leq 0$ and this contradicts the assumption because $\Delta_{m}(X,\mathcal{L})=\Delta_{m}(M,\mathcal{A})$ by \cite[Corollary 2.11]{Fukuma05}.
Therefore $h^{m}(\mathcal{O}_{M})>0$.
Hence by the above argument we get the assertion (2). 
Therefore we get the assertion of Claim \ref{CL0}. $\Box$
\\
\par
Since $\dim Y=3$ and we assume that $2g_{1}(X,\mathcal{L})-2=2\mathcal{L}^{n}$,
we have $\Delta_{3}(X,\mathcal{L})\leq 0$ by Claim \ref{CL0} if $(M,\mathcal{A})$ is a scroll over a normal projective variety of dimension $3$.
\par
Therefore we have $\Delta_{3}(X,\mathcal{L})\leq 0$ if $(X,\mathcal{L})$ satisfies (1) in Proposition \ref{PP2}.
\\
\\
(ii) Next we assume that $(X,\mathcal{L})$ satisfies (2.1) or (2.2) in Proposition \ref{PP2}.
\\
Assume that $(W,\mathcal{K})$ is the type 3 in \cite[Theorem 7.7.5]{BeSoBook}.
Then $2g_{1}(X,\mathcal{L})-2=2\mathcal{L}^{5}-1\neq 2\mathcal{L}^{5}$.
So we may assume that $(W,\mathcal{K})$ is not of this type.
\par
By Proposition \ref{PR3.1} and the assumption, we have $h^{j}(\mathcal{O}_{X})=h^{j}(\mathcal{O}_{M})=0$ and $g_{j}(X,\mathcal{L})=g_{j}(M,\mathcal{A})=0$ for every integer $j$ with $j\geq 3$.
By Remark \ref{Remark1.3} (2) we have
$\Delta_{n}(X,\mathcal{L})=\Delta_{n}(M,\mathcal{A})=h^{n}(\mathcal{O}_{M})-h^{n}(\mathcal{A})=-h^{n}(\mathcal{A})$.
Moreover by Remark \ref{Remark1.3} (3) and (4) we have
$$\Delta_{j}(X,\mathcal{L})=\Delta_{j}(M,\mathcal{A})
=g_{j}(M,\mathcal{A})-\Delta_{j+1}(M,\mathcal{A})+(n-j)h^{j}(\mathcal{O}_{M})-h^{j}(\mathcal{A}).
$$
So in order to calculate $\Delta_{3}(X,\mathcal{L})$, we have to calculate $h^{j}(\mathcal{A})$ with $j\geq 3$.
Then by Lemma \ref{LM2} we see that $h^{j}(\mathcal{A})=h^{j}(\mathcal{D})$ for every $j\geq 3$.
\\
(ii.1) If $(W,\mathcal{K})\cong (\mathbb{P}^{6},\mathcal{O}_{\mathbb{P}^{6}}(1))$,
then $\mathcal{D}=\mathcal{O}_{\mathbb{P}^{6}}(2)$ and we see that $h^{j}(\mathcal{D})=0$ for every $j\geq 2$.
\\
(ii.2) Assume that $(W,\mathcal{K})$ is the type 1 in \cite[Theorem 7.7.5]{BeSoBook}, that is, $(W,\mathcal{K})\cong (\mathbb{Q}^{5},\mathcal{O}_{\mathbb{Q}^{5}}(1))$.
Then $K_{W}$ is a Cartier divisor and $\mathcal{K}=K_{W}+3\mathcal{D}$.
Hence $3\mathcal{D}$ is also Cartier.
On the other hand $2\mathcal{D}$ is Cartier by \cite[Lemma 7.5.8]{BeSoBook}.
Hence $\mathcal{D}=3\mathcal{D}-2\mathcal{D}$ is also Cartier and $\mathcal{D}=\mathcal{O}_{\mathbb{Q}^{5}}(2)$.
Therefore by the Kawamata-Viehweg vanishing theorem \cite[Theorem 1-2-5]{KaMaMa}, we have $h^{j}(\mathcal{D})=h^{j}(\mathcal{O}_{\mathbb{Q}^{5}}(2))=h^{j}(K_{W}+\mathcal{O}_{\mathbb{Q}^{5}}(7))=0$ for every $j\geq 1$.
\\
(ii.3) Assume that $(W,\mathcal{K})$ is the type 2 in \cite[Theorem 7.7.5]{BeSoBook}.
Let $\pi:W\to C$ be the $\mathbb{P}^{4}$-bundle over a smooth curve $C$.
Then by \cite[Proposition 3.2.1]{BeSoBook}, there exists an ample vector bundle $\mathcal{E}$ on $C$ such that $W\cong \mathbb{P}_{C}(\mathcal{E})$ and $H(\mathcal{E})=\mathcal{K}$.
Since $W$ is smooth in this case, $\mathcal{D}$ is a Cartier divisor
and $\mathcal{D}=2H(\mathcal{E})+\pi^{*}(\mathcal{B})$ for $\mathcal{B}\in \mbox{Pic}(C)$.
\begin{s-Claim}\label{CL1}
$h^{j}(2H(\mathcal{E})+\pi^{*}(\mathcal{B}))=0$ for every $j\geq 2$.
\end{s-Claim}
\noindent{\em Proof.}
Since $R^{k}\pi_{*}(2H(\mathcal{E})+\pi^{*}(\mathcal{B}))=R^{k}\pi_{*}(2H(\mathcal{E}))\otimes\mathcal{B}=0$ for every positive integer $k$,
we see that $h^{j}(2H(\mathcal{E})+\pi^{*}(\mathcal{B}))=h^{j}(\pi_{*}(2H(\mathcal{E})+\pi^{*}(\mathcal{B})))$ for every $j\geq 0$.
Since $\dim C=1$, we get $h^{j}(\pi_{*}(2H(\mathcal{E})+\pi^{*}(\mathcal{B})))=0$ for every integer $j$ with $j\geq 2$.
Hence we get the assertion of Claim \ref{CL1}. $\Box$
\\
\par
By the above argument 
we have $h^{j}(\mathcal{A})=0$ for every integer $j$ with $j\geq 3$.
Hence by Remark \ref{Remark1.3} (4) we see that $\Delta_{j}(X,\mathcal{L})=\Delta_{j}(M,\mathcal{A})=0$ for every $j\geq 3$.
Therefore $\Delta_{3}(X,\mathcal{L})=0$
if $(X,\mathcal{L})$ satisfies (2.1) or (2.2) in Proposition \ref{PP2}.
\\
\par
From (i) and (ii) above, we see that $\Delta_{3}(X,\mathcal{L})\leq 0$ if $\kappa(K_{X}+(n-3)\mathcal{L})=-\infty$ and $2g_{1}(X,\mathcal{L})-2=2\mathcal{L}^{n}$.
Therefore we get the assertion of Lemma \ref{LM3}. $\Box$
\\
\\
From (III.1) and (III.2) we see that $C(3,3)$ implies $C(3,1)$.
\\
\\
(IV) Assume that $C(3,5)$.
Then we will prove that $C(3,1)$ holds.
First we see that $\kappa(K_{X}+(n-3)\mathcal{L})\leq 0$ by the same reason 
of the case (III) above.
\\
(IV.1) Assume that $\kappa(K_{X}+(n-3)\mathcal{L})=0$.
Then by the same argument as (III.1) above, we get $K_{X}+(n-3)\mathcal{L}=\mathcal{O}_{X}$
in this case.
\\
(IV.2) Assume that $\kappa(K_{X}+(n-3)\mathcal{L})=-\infty$.
By Proposition \ref{PR3.1} we see that $g_{3}(X,\mathcal{L})=0$ unless 
$(M,\mathcal{A})$ is a scroll over a normal projective variety of dimension $3$.
Next we assume that $(M,\mathcal{A})$ is a scroll over a normal projective variety $Y$ of dimension $3$.
Here we note that $Y$ is smooth and $(M,\mathcal{A})$ is a classical scroll over $Y$ by Proposition \ref{PP4}.
Then $g_{3}(X,\mathcal{L})=g_{3}(M,\mathcal{A})=h^{3}(\mathcal{O}_{M})=h^{3}(\mathcal{O}_{X})$ (see Remark \ref{TR2-3.1}).
If $h^{3}(\mathcal{O}_{X})>0$, then by Claim \ref{CL0} we have $2g_{1}(X,\mathcal{L})-2>2\mathcal{L}^{n}$ and this contradicts the assumption $C(3,5)$.
Hence $h^{3}(\mathcal{O}_{X})=0$, and $g_{3}(X,\mathcal{L})=0$ also holds in this case.
Therefore there does not any $(X,\mathcal{L})$ with $\kappa(K_{X}+(n-3)\mathcal{L})=-\infty$, $g_{3}(X,\mathcal{L})>0$ and $2g_{1}(X,\mathcal{L})-2=2\mathcal{L}^{n}$.
\\
\par
By (IV.1) and (IV.2) we get $C(3,1)$. Therefore these complete the proof of Theorem \ref{MT3}.
$\Box$
\\
\par
By Theorems \ref{MT4}, \ref{MT2} and \ref{MT3} we get the following corollary.

\begin{Corollary}\label{CR3}
Let $(X,\mathcal{L})$ be a polarized manifold of dimension $n\geq 5$.
Assume that $\dim \mbox{\rm Bs}|\mathcal{L}|=2$. 
Then Conjecture {\rm \ref{CJ}} is true.
\end{Corollary}

\begin{s-Remark}
Next we consider the case where $n=4$.
By the same argument as the proof of Theorem \ref{MT3}, we can check that the following implications hold:
$C(3.1)\Longrightarrow C(3.2)$, $C(3.1)\Longrightarrow C(3.4)$,
$C(3.2)\Longrightarrow C(3.4)$, and 
$C(3.3)\Longrightarrow C(3.5)$.
\par
Next we consider the implication $C(3.3)\Longrightarrow C(3.1)$.
By the same argument as in the case (III) in Theorem \ref{MT3} we have $\kappa(K_{X}+\mathcal{L})\leq 0$.
\par
If $\kappa(K_{X}+\mathcal{L})=0$, then we can prove that $K_{X}+\mathcal{L}=\mathcal{O}_{X}$.
So we assume that $\kappa(K_{X}+\mathcal{L})=-\infty$.
Since $n=4$, we see that $(X,\mathcal{L})$ satisfies (2) in Theoem \ref{TH2}.
Then $h^{4}(\mathcal{O}_{M})=0$ and $h^{4}(\mathcal{A})=0$ because $h^{4}(\mathcal{A})=h^{0}(K_{M}-\mathcal{A})$ and $\kappa(M)=-\infty$.
By Proposition \ref{PR3.1} and Lemma \ref{LM1} we have $g_{3}(X,\mathcal{L})=g_{3}(M,\mathcal{A})=0$ and $h^{3}(\mathcal{O}_{X})=h^{3}(\mathcal{O}_{M})=0$ unless
$(M,\mathcal{A})$ is a scroll over a normal projective variety of dimension $3$.
Therefore $\Delta_{4}(M,\mathcal{A})=0$ and $\Delta_{3}(M,\mathcal{A})=-h^{3}(\mathcal{A})\leq 0$ unless $(M,\mathcal{A})$ is a scroll over a normal $3$-fold.
Since $\Delta_{3}(X,\mathcal{L})=\Delta_{3}(M,\mathcal{A})$, we get $\Delta_{3}(X,\mathcal{L})\leq 0$.
But this contradicts the assumption.
Therefore if $(M,\mathcal{A})$ is not a scroll over a normal $3$-fold,
then $C(3.3)$ implies $C(3.1)$.
\par
If $(M,\mathcal{A})$ is a classical scroll over a smooth $3$-fold,
then by Claim \ref{CL0} we also see that $C(3.3)$ implies $C(3.1)$.
By the same argument as above, we see that $C(3.5)$ implies $C(3.1)$
if $(M,\mathcal{A})$ is a classical scroll over a smooth $3$-fold.
\par
So in order to prove that Conjecture \ref{CJ} for $n=4$ and $i=3$ is true, it suffices to consider the case where a reduction of $(X,\mathcal{L})$ is a scroll over a normal projective variety $Y$ with $\dim Y=3$, but not a classical scroll over a smooth $3$-fold $Y$.
\end{s-Remark}

\subsection{Some remarks}

Finally we would like to give a comment about Conjecture \ref{CJ}.
\\
\\
(a)
We can easily see that $C(i,1)$ implies $C(i,2)$ and $C(i,4)$, and $C(i,2)$ (resp. $C(i,4)$) implies $C(i,3)$ (resp. $C(i,5)$) by the same argument as the proof of Theorem \ref{MT4}.
\\
\\
(b) By looking at the proof of Theorems \ref{MT4}, \ref{MT2} or \ref{MT3}, we can prove the following:

\begin{s-Proposition}\label{PP1}
If there does not exist any polarized manifold $(X,\mathcal{L})$ with $\kappa(K_{X}+(n-i)\mathcal{L})=-\infty$, $\Delta_{i}(X,\mathcal{L})>0$ {\rm (}resp. $g_{i}(X,\mathcal{L})>0${\rm )} and $2g_{1}(X,\mathcal{L})-2=(i-1)\mathcal{L}^{n}$, then we see that $C(i,3)$ {\rm (}resp. $C(i,5)${\rm )} implies $C(i,1)$.
\end{s-Proposition}
So it is important to know whether there exists an example of $(X,\mathcal{L})$ with $\kappa(K_{X}+(n-i)\mathcal{L})=-\infty$, $\Delta_{i}(X,\mathcal{L})>0$ (resp. $g_{i}(X,\mathcal{L})>0$) and $2g_{1}(X,\mathcal{L})-2=(i-1)\mathcal{L}^{n}$.
\\
\\
(c) We can regard the following result as the case where $i=n$ in Conjecture \ref{CJ}.
\begin{s-Proposition}\label{PP3}
Let $X$ be a smooth projective variety of dimension $n$.
Then the following are equivalent one another.
\begin{itemize}
\item [\rm (i)] $K_{X}\sim \mathcal{O}_{X}$.
\item [\rm (ii)] $\Delta_{n}(X,\mathcal{L})=1$ and $2g(X,\mathcal{L})-2=(n-1)\mathcal{L}^{n}$ hold for any ample line bundle $\mathcal{L}$.
\item [\rm (iii)] $\Delta_{n}(X,\mathcal{L})>0$ and $2g(X,\mathcal{L})-2=(n-1)\mathcal{L}^{n}$ hold for any ample line bundle $\mathcal{L}$.
\item [\rm (iv)] $g_{n}(X,\mathcal{L})=1$ and $2g(X,\mathcal{L})-2=(n-1)\mathcal{L}^{n}$ hold  for any ample line bundle $\mathcal{L}$.
\item [\rm (v)] $g_{n}(X,\mathcal{L})>0$ and $2g(X,\mathcal{L})-2=(n-1)\mathcal{L}^{n}$ hold  for any ample line bundle $\mathcal{L}$.
\end{itemize}
\end{s-Proposition}
\noindent{\em Proof.}
(i)$\Rightarrow$(ii): By Remark \ref{Remark1.3} (2) we have $\Delta_{n}(X,\mathcal{L})=h^{n}(\mathcal{O}_{X})-h^{n}(\mathcal{L})$.
By assumption we get $h^{n}(\mathcal{O}_{X})=h^{0}(K_{X})=1$.
Next we calculate $h^{n}(\mathcal{L})$.
Since $h^{n}(\mathcal{L})=h^{0}(K_{X}-\mathcal{L})=h^{0}(-\mathcal{L})$, we see that $h^{n}(\mathcal{L})=0$ because $\mathcal{L}$ is ample.
Therefore $\Delta_{n}(X,\mathcal{L})=1$.
Of course $2g(X,\mathcal{L})-2=(n-1)\mathcal{L}^{n}$ holds since $K_{X}\sim \mathcal{O}_{X}$.
\\
(ii)$\Rightarrow$(iii): This is trivial.
\\
(iii)$\Rightarrow$(i): 
Since $0<\Delta_{n}(X,\mathcal{L})=h^{n}(\mathcal{O}_{X})-h^{n}(\mathcal{L})$, we have $h^{n}(\mathcal{O}_{X})\geq 1$.
By the Serre duality we see that $h^{0}(K_{X})\geq 1$.
On the other hand, since $2g(X,\mathcal{L})-2=(n-1)\mathcal{L}^{n}$, we have $K_{X}\mathcal{L}^{n-1}=0$.
Therefore we get $K_{X}\sim \mathcal{O}_{X}$.
\\
(i)$\Rightarrow$(vi): Since $K_{X}\sim \mathcal{O}_{X}$, we see that $2g(X,\mathcal{L})-2=(K_{X}+(n-1)\mathcal{L})\mathcal{L}^{n-1}=(n-1)\mathcal{L}^{n}$ and $h^{0}(K_{X})=1$.
By the Serre duality we have $h^{n}(\mathcal{O}_{X})=h^{0}(K_{X})=1$ holds.
Hence $g_{n}(X,\mathcal{L})=h^{n}(\mathcal{O}_{X})=1$ by Remark \ref{Remark1} (2).
\\
(vi)$\Rightarrow$(v): This is trivial.
\\
(v)$\Rightarrow$(i): By assumption, we have $h^{n}(\mathcal{O}_{X})=g_{n}(X,\mathcal{L})>0$.
We also note that $K_{X}\mathcal{L}^{n-1}=0$ by the assumption that $2g(X,\mathcal{L})-2=(n-1)\mathcal{L}^{n}$. Hence we see that $K_{X}\sim \mathcal{O}_{X}$ because $h^{0}(K_{X})=h^{n}(\mathcal{O}_{X})>0$. $\Box$

\begin{flushright}
Department of Mathematics \\
Faculty of Science \\
Kochi University \\
Akebono-cho, Kochi 780-8520 \\
Japan \\
E-mail: fukuma@kochi-u.ac.jp \\
\end{flushright}
\end{document}